\documentclass[11pt,reqno]{amsart}

\usepackage[text={160mm,240mm},centering]{geometry}            
\usepackage{graphicx}
\usepackage{listings}
\usepackage{amssymb}
\usepackage{dsfont}
\usepackage{url}
\usepackage{verbatim}
\usepackage{bbm}

\usepackage{amssymb,amsmath,setspace,geometry,indentfirst,changepage,inputenc}
\usepackage{mathrsfs,amsfonts,savesym,graphicx,bm,amsthm}

\usepackage{xurl}

\usepackage[all]{xy}
\usepackage{mathrsfs}
\usepackage{hyperref}
\usepackage{color}

\usepackage{dcolumn}
\usepackage{pdfpages}

\usepackage{booktabs}
\usepackage{diagbox} 
\usepackage{multirow} 
\usepackage{makecell} 
\usepackage{longtable} 
\usepackage{array}
\usepackage{tabularx}
\usepackage{color}

\usepackage{amsmath,amsthm}

\usepackage[mathscr]{eucal}
\usepackage[all]{xy}
\usepackage{mathrsfs}
\usepackage{hyperref}
\usepackage{color}
\usepackage{fancyhdr}

\newtheorem{theorem}{Theorem}[section]
\newtheorem{lemma}[theorem]{Lemma}

\newtheorem{conjecture}[theorem]{Conjecture}
\newtheorem{theorem-definition}[theorem]{Theorem-Definition}
\newtheorem{corollary}[theorem]{Corollary}

\newtheorem{proposition}[theorem]{Proposition}

\newtheorem{question}[theorem]{Question}
\newtheorem{definition-theorem}[theorem]{Definition-Theorem}

\newtheorem{theorem-defintion}[theorem]{Theorem-Definition}
\newtheorem{corollary-definition}[theorem]{Corollary-Definition}

\theoremstyle{definition}
\newtheorem{example}[theorem]{Example}
\newtheorem{definition}[theorem]{Definition}

\newtheorem{notation}[theorem]{Notation}

\newtheorem{remark}[theorem]{Remark}

\newtheorem*{ack}{Acknowledgements}

\author{Quan Shi}
\address{Department of Mathematical Sciences, Tsinghua University, Beijing, 100084, P. R. China.}
\email{shiq24@mails.tsinghua.edu.cn / thusq20@gmail.com}

\author{Huaiqing Zuo}
\address{Department of Mathematical Sciences, Tsinghua University, Beijing, 100084, P. R. China.}
\email{hqzuo@mail.tsinghua.edu.cn}

\allowdisplaybreaks

\title{Variation of Archimedean Zeta Function and $n/d$-Conjecture for Generic Multiplicities}

\begin{document}

	\begin{abstract}
		For $f_1,...,f_r\in \mathbb C[z_1,...,z_n]\setminus \mathbb C$, we introduce the variation of archimedean zeta function. As an application, we show that the $n/d$-conjecture, proposed by Budur, Musta\c{t}\u{a}, and Teitler, holds for generic multiplicities. Consequently, strong monodromy conjecture holds for hyperplane arrangements with generic multiplicities as well.
		\vspace{0.5em}
		
		
	\end{abstract}
	
	\maketitle

	\tableofcontents

	\section{Introduction}

	Let $f \in \mathbb C[z_1,...,z_n]\setminus \mathbb C$. Two profound objects associated with $f$ can be constructed. The first is the \textit{$b$-function $b_f(s)$}, also called the Bernstein-Sato polynomial. Let $A_n$ be the Weyl algebra, i.e. the subring of $\mathrm{End}_{\mathbb C\text{-linear}}\big(\mathbb C[z_1,...,z_n]\big)$ generated by $\mathbb C, z_1,...,z_n,{\frac{\partial}{\partial z_1}},...,{\frac{\partial}{\partial z_n}}$. $s$ is a formal symbol, $A_n[s] := A_n\otimes_{\mathbb C} \mathbb C[s]$, and $f^{s}$ is a formal power of $f$.
	\begin{definition}\label{Definition of $b$-function}
		The (algebraic) $b$-function of $f$, denoted by $b_f(s)\in \mathbb C[s]$, is the monic polynomial of the smallest degree that admits a $P \in A_n[s]$ such that
		\begin{displaymath}
			P\cdot \big( f\cdot f^s \big) = b_f(s)\cdot f^s.
		\end{displaymath}
		Here $A_n[s]$ acts on $\mathbb C[z_1,...,z_n,\frac{1}{f},s]\cdot f^{s}$ in the expected way.
	\end{definition}
	
	The roots of $b_f(s)$ alone are significant in algebraic geometry. In \cite{Rationality_of_Roots_of_B-Function_Kashiwara}, Kashiwara proved that all roots of $b_f(s)$ are negative rational numbers. Later, Kashiwara and Malgrange independently proved that the exponents of the roots coincide with monodromy eigenvalues, see \cite{Kashiwara_Malgrange_Theorem_Kashiwara}, \cite{Kashiwara_Malgrange_Theorem_Malgrange} and see \cite{Milnor_Fibration_Monodromy_2}, \cite{Milnor_Fibration_Monodromy_3}, \cite{Mixed_Hodge_Applied_to_Singularity} for monodromy. 
	
	The second object is the \textit{motivic zeta function $Z_{f}^{\mathrm{mot}}(s)$}. It can be intrinsically defined by using arcs and jets and applying the motivic integration, see \cite{Motivic_Integration_on_Arbitrary_Varieties_Denef_Loeser}. Here, we leave the extrinsic construction provided in \cite{Motivic_Igusa_Zeta_Function}.
	\begin{theorem}\label{Denef-Loeser formula}
		For $f\in \mathbb C[z_1,...,z_n]\setminus \mathbb C$, let $\mu : Y \to {\mathbb C}^n$ be a log resolution of $f$. Suppose $\mu^*(\mathrm{div}\, f) = \sum_{i\in S} N_i D_i$ and $K_\mu = \sum_{i\in S}(\nu_i-1)D_i$, where $\{D_i\}_{i\in S}$ is a collection of simple normal crossings divisors and $K_\mu$ is the relative canonical divisor. $\mathbb L^{-s}$ is a formal power. Then the \textit{motivic zeta function of $f$} is given by
		\begin{displaymath}	
			Z_{f}^{\mathrm{mot}}(s) := \mathbb L^{-n} \sum_{I\subseteq S} [D_I^\circ] \cdot \prod_{i\in I}\frac{(\mathbb L-1)\mathbb L^{-(N_is+\nu_i)}}{1-\mathbb L^{-(N_is+\nu_i)}}. \tag{\ref{Denef-Loeser formula}.1}
		\end{displaymath}
		Here ${D_I^\circ = (\bigcap_{i\in I} D_i) \setminus (\bigcup_{j\in S\setminus I} D_j)}$, $[-]$ denotes the class of a variety in the Grothendieck ring $K_0[\mathrm{Var}_{\mathbb C}]$ of complex algebraic varieties, and $\mathbb L = [\mathbb A^1_{\mathbb C}]$. 
	\end{theorem}
	\begin{remark}
		The independence of resolutions of the expression can be shown by using the weak factorization theorem, see \cite{Torification_and_Factorization_of_Birational_Maps_Abramovich_Karu_Matsuki}. 
	\end{remark}
	From (\ref{Denef-Loeser formula}.1), one can define the poles of $Z_{f}^{\mathrm{mot}}(s)$ in a natural way, see \cite{Introduction_to_the_Monodromy_Conjecture_Veys_2024}.
	
	The following is the well-known \textit{monodromy conjecture}, see \cite{Motivic_Igusa_Zeta_Function}, which is one of the main open problems in singularity theory.
	\begin{conjecture}[monodromy conjecture]\label{Monodromy_Conjecture}
		For $f\in \mathbb C[z_1,...,z_n]\setminus \mathbb C$, if $s_0$ is a pole of $Z_{f}^{\mathrm{mot}}(s)$, then the following hold.
		
		\noindent(a) (Weak) $e^{2\pi \sqrt{-1}\cdot s_0}$ is a monodromy eigenvalue at some point of $\{f = 0\}$.
		
		\noindent(b) (Strong) $s_0$ is a root of $b_f(s)$.	
	\end{conjecture}
	Monodromy conjecture would bridge motivic objects with objects from geometry and $\mathscr D$-module theory, thus being very interesting. However, monodromy conjecture is very difficult even for concrete classes of polynomials. In this paper, an important class of polynomials, hyperplane arrangements, is concerned.
	
	A \textit{hyperplane arrangement} in $\mathbb C^n$ is the union of a finite collection of hyperplanes. Let $A = \{f = 0\}$ be a hyperplane arrangement in $\mathbb C^n$ with $f \in \mathbb C[z_1,...,z_n]\setminus \mathbb C$ and $A = \bigcup_{i=1}^r V_i$ with $V_i$ being hyperplanes. $A$ is called \textit{decomposable} if, after a linear coordinate change, $f$ can be written as the product of two non-constant polynomials in two disjoint sets of variables. Otherwise, $A$ is called \textit{indecomposable}. $A$ is called \textit{essential} if $\bigcap_{i = 1}^r V_i = \bm 0$. Otherwise, $A$ is called \textit{non-essential}. $A$ is called \textit{central} if $\bm 0\in V_i$ for all $i=1,...,r$. 
	
	
	In terms of the monodromy conjecture for hyperplane arrangements,  we have the following theorem and conjecture from \cite{Monodromy_Conjecture_for_Hyperplane_Arrangement_Budur_Mustata_Teitler}.
	\begin{theorem}\label{Weak Monodromy Conjecture for hyperplane arrangements} 
		Weak monodromy conjecture holds for hyperplane arrangements. Moreover, if the $n/d$-conjecture, i.e. {Conjecture \ref{Proof_of_Weak_MC_hyperplane_arrangement}} holds, then the strong monodromy conjecture holds for hyperplane arrangements.
	\end{theorem}
	\begin{conjecture}[$n/d$-conjecture]\label{Proof_of_Weak_MC_hyperplane_arrangement}
		Let $f \in \mathbb C[z_1,...,z_n]\setminus \mathbb C$ be homogeneous of degree $d$ and suppose that $\{f = 0\}$ defines a central essential indecomposable hyperplane arrangement, then $b_f(-\frac{n}{d}) = 0$.
	\end{conjecture}
	The $n/d$-conjecture is important not only for its application to the strong monodromy conjecture but also for its combinatoric nature. In general, the $b$-function is not a combinatoric invariant for hyperplane arrangements, see \cite{Jacobian module}, but the conjecture predicts that it has a root $-n/d$, which can be read from the combinatorics.
	
	\vspace{0.5em}

	In this paper, we develop a new object, \textit{variation of archimedean zeta function}, in Subsection \ref{Subsection_Variation_of_Archimedean_Zeta_Function}, which is a suitable tool to study the $b$-function in the multiplicative family $\{f_1^{a_1}\cdots f_r^{a_r} \mid  (a_1,...,a_r) \in \mathbb N_{>0}^r\}$, where $f_1,...,f_r\in \mathbb C[z_1,...,z_n]\setminus\mathbb C$ are fixed polynomials. Using the variation of archimedean zeta function, we obtain the following result about $n/d$-conjecture.
	
	\begin{theorem}\label{Generically_n_over_d}
		Let $f_1,...,f_r\in \mathbb C[z_1,...,z_n]\setminus \mathbb C$ be homogeneous polynomials of degree one. Suppose $A = \{f_1f_2\cdots f_r = 0\}$ is a central essential indecomposable hyperplane arrangement.  Let $\varphi\in C_c^{\infty}(\mathbb C^n) $ be such that $\varphi(\bm 0) > 0$ and $\varphi \geq 0$. Then there is a non-empty analytically open subset $U \subseteq \{(a_1,...,a_r) \in \mathbb C^r \mid \sum_{i=1}^r a_i = 1\}$ such that for $\bm a = (a_1,...,a_r) \in \mathbb N_{>0}^r$, if $(\frac{a_1}{\sum_{i=1}^r a_i},...,\frac{a_r}{\sum_{i=1}^r a_i}) \in U$, then $-n/\sum_{i=1}^r a_i$ is a pole of $Z_{\bm f^{\bm a},\varphi}(s)$, where $\bm f^{\bm a} := f_1^{a_1}\cdots f_r^{a_r}$ and $Z_{\bm f^{\bm a},\varphi}(s)$ is defined after this theorem. In particular, for such $\bm a$, $-n/\sum_{i=1}^r a_i$ is a root of $b_{\bm f^{\bm a}}(s)$. That is, the $n/d$-conjecture holds for generic multiplicities.
	\end{theorem}
	
	%

	The $Z_{\bm f^{\bm a},\varphi}(s)$ in the above theorem is the \textit{archimedean zeta function} with respect to $\bm f^{\bm a}$ and $\varphi$. That is, for $f\in \mathbb C[z_1,...,z_n]\setminus \mathbb C$ and $\varphi \in C_c^{\infty}(\mathbb C^n)= C_c^{\infty}(\mathbb C^n,\mathbb C)$, $Z_{f,\varphi}(s)$ is defined as the meromorphic extension of the following function of $s$. Here, we adopt a normalization by $(\frac{\sqrt{-1}}{2})^n$.
	\begin{displaymath}
		\mathbb U \to \mathbb C,\ s\mapsto (\frac{\sqrt{-1}}{2})^n \int_{\mathbb C^n} \vert f\vert^{2s} \varphi\, \mathrm{d}z_1\wedge \mathrm{d} \bar z_1\wedge \cdots \wedge \mathrm{d}z_{n} \wedge \mathrm{d} \bar z_{n},
	\end{displaymath}
	where $\mathbb U := \{s \in \mathbb C \mid \mathrm{Re}\, s > 0\}$. The existence of such meromorphic extension was settled by Gel'fand and Bernstein, and Atiyah independently, see \cite{Continuation_Problem_of_Zeta_Function_Bernstein_and_Gelfand}, \cite{Continuation_Problem_of_Zeta_Function_Atiyah}, \cite{Bernstein-polynomial_Continuations}. Log resolutions were applied in \cite{Continuation_Problem_of_Zeta_Function_Bernstein_and_Gelfand}, \cite{Continuation_Problem_of_Zeta_Function_Atiyah}, while \cite{Bernstein-polynomial_Continuations} used $b$-functions. From the proof in \cite{Bernstein-polynomial_Continuations}, we know that if $\alpha \in [-1,0)$ is a pole of $Z_{f,\varphi}(s)$, then $\alpha$ is a root of $b_f(s)$. It is a classical strategy in singularity theory. We put more related information in Subsection \ref{Subsection_Archimedean_Zeta_Function}. 
	
	\vspace{0.5em}

	As a corollary of {Theorem \ref{Generically_n_over_d}}, a generic result about the strong monodromy conjecture for hyperplane arrangements is obtained.

	\begin{corollary}\label{generic monodromy conjecture}
		Let $f_1,...,f_r\in \mathbb C[z_1,...,z_n]\setminus \mathbb C$ be degree-one polynomials. Suppose $A = \{f_1f_2\cdots f_r = 0\}$ is a hyperplane arrangement, not necessarily central. Then there exists an analytically open subset $U \subseteq \{(a_1,...,a_r) \in \mathbb U^r \mid \sum_{i=1}^r a_i = 1\}$, such that if $\bm a = (a_1,...,a_r) \in \mathbb N_{>0}^r$ and $(\frac{a_1}{\sum_{i=1}^r a_i},...,\frac{a_r}{\sum_{i=1}^r a_i}) \in U$, then the strong monodromy conjecture holds for $\bm f^{\bm a} = f_1^{a_1}\cdots f_r^{a_r}$. That is, the strong monodromy conjecture for hyperplane arrangements holds for generic multiplicities.
	\end{corollary}
	
	%

	In fact, Theorem \ref{Generically_n_over_d} is a consequence of Theorem \ref{generic existence} below, which is a more general result about poles of archimedean zeta functions.  Let $f_1,...,f_r \in \mathbb C[z_1,...,z_n]\setminus \mathbb C$ and $\mu : Y \to \mathbb C^n$ be a log resolution of $\prod_{i=1}^rf_i$. Suppose $\{E_j\}_{j\in I}$ is a collection of simple normal crossings divisors on $Y$ such that $\mu^{-1}(\{f_1\cdots f_r = 0\}) = \bigcup_{j\in I} E_j$, $\mu^*(\mathrm{div}\, f_i) = \sum_{j\in I} N_{ij} E_j$, and $K_\mu = \sum_{j\in I} (\nu_j - 1)E_j$. Let $P_j : \mathbb C^r \to \mathbb C,\ (a_1,...,a_r) \mapsto \sum_{i=1}^r N_{ij}a_i$. A tuple $\bm a_0 = (a_{01},...,a_{0r}) \in \mathbb N_{>0}^r$ is called a \textit{good tuple} at $j\in I$ if $-\frac{P_{j'}(\bm a_0)}{P_{j}(\bm a_0)}\nu_j+\nu_{j'} > 0$ for all $j'\neq j$. In particular, $\frac{\nu_j}{P_j(\bm a_0)}$ is the log canonical threshold of $\bm f^{\bm a_0} = f_1^{a_{01}}\cdots f_r^{a_{0r}}$, see Remark \ref{lct_define} for definition.
	
	\begin{theorem}\label{generic existence}
		The notation is as above. Let $\varphi \in C_c^{\infty}(\mathbb C^n)$ with $\varphi \geq 0$. Suppose $\bm a_0 \in \mathbb N_{>0}^r$ is a good tuple at $j \in I$ and $\varphi|_{\mu(E_j)} \neq 0$. Then there is a nonempty analytically open subset $U \subseteq \{\bm b\in \mathbb C^r \mid P_j(\bm b) = 1\} \simeq \mathbb C^{r-1}$ such that for $\bm a = (a_1,...,a_r) \in \mathbb N_{>0}^r$, if $(\frac{a_1}{P_j(\bm a)},...,\frac{a_r}{P_j(\bm a)}) \in U$, then $-\frac{\nu_j}{P_j(\bm a)}$ is a pole of $Z_{\bm f^{\bm a},\varphi}(s)$, where $\bm f^{\bm a} := f_1^{a_1}\cdots f_r^{a_r}$.  In particular, for such $\bm a$, if $-\frac{\nu_j}{P_j(\bm a)} \geq -1$, then $-\frac{\nu_j}{P_j(\bm a)}$ is a root of $b_{\bm f^{\bm a}}(s)$.
	\end{theorem}
	
	%

	We can enhance the result of Theorem \ref{Generically_n_over_d} when $n = 2$. After a deeper analysis on the residue of the archimedean zeta function at $-2/d$, a stronger outcome is achieved. 
	
	
	\begin{theorem}\label{n_over_d_in_Dimension_two}
		Let $f_1,...,f_r\in \mathbb C[z_1,z_2]\setminus \mathbb C$ be degree-one polynomials. Suppose $A = \{f_1f_2\cdots f_r = 0\}$ is a central essential indecomposable hyperplane arrangement.  Let $\alpha,\beta\in C_c^{\infty}(\mathbb C)$ be such that $\alpha,\beta \geq 0$ and there exists a neighborhood $D$ of $\,0\in \mathbb C$, $\alpha|_D = \beta|_D \equiv 1$. Take $\varphi \in C_c^{\infty}(\mathbb C^2)$ such that $\varphi(z_1,z_2) := \alpha(z_1)\beta(z_2)$ for all $(z_1,z_2) \in \mathbb C^2$.
		
		Then for all $\bm a = (a_1,...,a_r) \in \mathbb N_{>0}^r$, the order of $-2/d$ as a pole of $Z_{\bm f^{\bm a},\varphi}(s)$ is not greater than $2$, where $\bm f^{\bm a} := f_1^{a_1}\cdots f_r^{a_r}$. Moreover, if $-2/d$ is not an order two pole of $Z_{\bm f^{\bm a},\varphi}(s)$, then it is an order one pole with a positive (negative resp.) residue if $\max \{a_i\}_{i=1}^r \leq d/2$ ($\max \{a_i\}_{i=1}^r > d/2$ resp.).
	\end{theorem}
	
	\begin{remark}
		Combing {Remark \ref{tame leq 3}} with  \cite[Theorem 1.9]{Quelques}, one can also deduce that $-2/d$ is a pole of $Z_{f,\varphi}(s)$ for some $\varphi$. However, the statement of Theorem \ref{n_over_d_in_Dimension_two} is stronger since a specific $\varphi$ is provided and we find the sign of the residue.
	\end{remark}
	
	For central essential indecomposable hyperplane arrangements, some existing results about $b$-functions present more behaviors of their roots. Thus, we propose the following question.
	\begin{question}[stronger $n/d$-conjecture]\label{Strong_n_over_d_Conjecture}
		Let $f \in \mathbb C[z_1,...,z_n]\setminus \mathbb C$ be homogeneous of degree $d$ and suppose that $\{f = 0\}$ defines a central essential indecomposable hyperplane arrangement. Then are $-\frac{n}{d},-\frac{n+1}{d},...,-\frac{d-1}{d}$ all roots of $b_f(s)$?
	\end{question}
	For example, according to the main theorems of \cite{B_S_for_Centeral_Generic_Hyperplane_Arrangements_Walther}, \cite{B_S_for_Centeral_Generic_Hyperplane_Arrangements_Saito}, and \cite{Tame_Arrangement}, this problem has a positive answer when $f$ is reduced and generic, or $f$ is tame. 
	
	\vspace{0.5em}

	The paper is organized as follows. Section \ref{section2} is for preliminary materials. In Section \ref{section3}, we develop the variation of archimedean zeta function and prove Theorem \ref{generic existence}, Theorem \ref{Generically_n_over_d}, Corollary \ref{generic monodromy conjecture} in this order. In Section \ref{section4}, we prove Theorem \ref{n_over_d_in_Dimension_two}. 
	
	\vspace{1em}

	The following is some notation adopted in the main body of this paper.
	\vspace{-0.5em}
	\begin{notation}
		\label{N&C} (1) We use a boldface letter to mean a tuple, with its size and the types of entries depending on context. (2) $\bm z$ denotes the complex coordinates $(z_1,...,z_n)$, and we set $\mathrm{d} \bm z := (\frac{\sqrt{-1}}{2})^n\mathrm{d}z_1\wedge \mathrm{d} \bar z_1\wedge \cdots \wedge \mathrm{d}z_{n} \wedge \mathrm{d} \bar z_{n}$. (3)  For polar coordinates $\{(\rho_i,\theta_i)\}_{i=1,...,n}$, we write $\mathrm{d} \rho_1\cdots\mathrm{d} \rho_n\mathrm{d} \theta_1\cdots\mathrm{d}\theta_n$ as $\mathrm{d} \bm \rho\mathrm{d} \bm \theta$. (4) For functions $f_1,...,f_r$ and $\bm a \in\mathbb N^r$, we write $\bm f^{\bm a} := \prod_{i=1}^r f_i^{a_i}$. (5) For an array $\bm m \in \mathbb Z^k$ for some $k$, $\vert m\vert$ denotes $\sum_{i=1}^k m_i$. (6) $\mathbb D_e(a)$ denotes the disk in $\mathbb C$ with center $a$ and radius $e$ and $\mathbb B(\bm x,\varepsilon)$ be the ball in a general metric space with center $\bm x$ and radius $\varepsilon$. (7) $\mathbb U := \{s \in \mathbb C \mid \mathrm{Re}\, s > 0\}$.
	\end{notation}
	
	\vspace{0.5em}

	
	\begin{ack}
		Zuo was supported by NSFC Grant 12271280. The authors thank Nero Budur for suggestions on this paper and his course in Spring 2024, during which they learned a lot more about the $n/d$-conjecture and the monodromy conjecture. They also thank the referee for very helpful advice.
	\end{ack}

	\section{Preliminary}\label{section2}
	
	\subsection{$b-$function}
	
	Here, we review several important properties of $b$-functions. Most of  the following can be found in the lecture notes of Popa \cite{DMBG}. Let $f \in \mathbb C[\bm z]\setminus \mathbb C = \mathbb C[z_1,...,z_n]\setminus \mathbb C$.

	\begin{theorem}[\cite{LCT_LARGEST_POLE}]\label{LCT_B_ROOT}
		The largest root of the Bernstein-Sato polynomial coincides with the negative of $\mathrm{lct}(f)$, the log canonical threshold of $f$.
	\end{theorem}
	\begin{remark}\label{lct_define}
		(1) For a subset $Z\subseteq \mathbb C^n$, $\mathrm{lct}_Z(f) \in \mathbb R_{>0}$ is the largest number $s$ such that $\vert f\vert^{-2s}$ is locally integrable near $Z$ and $\mathrm{lct}(f) = \mathrm{lct}_{\mathbb C^n}(f)$. $\mathrm{lct}(f)$ ($\mathrm{lct}_Z(f)$ resp.) is called the \textit{log canonical threshold} of $f$ (of $f$ at $Z$ resp.). (2) More generally, roots of $b_f(s)$ can be further related to jumping coefficients of multiplier ideals, see \cite{Jumping_Coefficients_of_Multiplier_Ideals}. 
	\end{remark}
	The following is the Kashiwara-Malgrange theorem, cf. \cite{Kashiwara_Malgrange_Theorem_Kashiwara} and \cite{Kashiwara_Malgrange_Theorem_Malgrange}. It relates $b$-functions with monodromies, thus providing a geometrical interpretation of the $b$-function.
	\begin{theorem}\label{Ch1_Sec2_Kashiwara_Malgrange_Theorem}
		Let $f\in \mathbb C[\bm z]\setminus \mathbb C$ be a non-constant polynomial. Suppose $\alpha$ is a root of $b_f(s)$. Then $e^{2\pi \sqrt{-1} \alpha}$ is a monodromy eigenvalue at some point of $\{f = 0\}$. Conversely, eigenvalues of every local monodromy are of this type.
	\end{theorem}
	\begin{remark}
		$b$-function can also be defined for analytic germs, only by substituting $\mathbb C[\bm z]$ and $A_n$ in {Definition \ref{Definition of $b$-function}} with the ring of convergent power series $\mathbb C\{\bm z\}$ and $\mathbb C\{\bm z\}[\bm{\frac{\partial}{\partial z}}]$. In this paper, if not emphasizing ``in analytic sense'', all $b$-functions are in algebraic sense.
	\end{remark}
	\begin{remark}
		In both algebraic and analytic senses, $-1$ is a root of $b_a(s)$ as long as $a$ is not a constant. Define the \textit{reduced $b$-function} $\tilde b_a(s) := b_a(s)/(s+1)$ and its set of roots $\tilde R_a := \{\tilde b_a = 0\} \subseteq \mathbb C$.
	\end{remark}
	\begin{remark}
		The general computation of $b$-functions is very difficult. However, this invariant admits some Thom-Sebastiani-type properties. Suppose $f\in \mathbb C[x_1,...,x_n]\setminus \mathbb C$ and $g\in \mathbb C[y_1,...,y_m]\setminus \mathbb C$. Then (1) $b_{fg}(s) = b_{f}(s)\cdot b_g(s)$, see \cite{shi2024tensorpropertybernsteinsatopolynomial} and \cite{lee2024multiplicativethomsebastianibernsteinsatopolynomials}, and (2) in analytic setting, we have $\tilde R_{f}+\tilde R_{g} \subseteq \tilde R_{f+g} + \mathbb Z_{\leq 0}$ and $\tilde R_{f+g} \subseteq \tilde R_{f} + \tilde R_{g} + \mathbb Z_{\geq 0}$, see \cite{microlocal b-fun}.
	\end{remark}
	
	\begin{remark}
		The local $b$-function is defined as follows. For a principal affine open subset $U = D(g) \subseteq \mathbb C^n$, one defines $b_{f|_U}(s)$ as the minimal monic polynomial that admits a solution $P \in A_n[\frac{1}{g},s]$ such that $P\cdot \big( f\cdot f^s\big) = b_{f|_U}(s)\cdot f^s$.
		For $x \in \mathbb C^n$, one defines the local $b$-function $b_{f,x}(s) := \mathrm{gcd}_{x \in U}\, b_{f|_U}(s)$, where $U$ runs through all principal affine open subsets containing $x$.
	\end{remark}
	
	\begin{remark}\label{basic facts about B-fun}
		Here are three basic facts.
		
		\noindent(1) If $U = D(g)$ and $h$ is an invertible element of $\mathbb C[\bm z,\frac{1}{g}]$, then $b_{f|_U}(s) = b_{h\cdot f|_U}(s)$.
		
		\noindent(2) $b_{f}(s) = \mathrm{lcm}_{x\in \{f = 0\}} b_{f,x}(s)$.
		
		\noindent(3) If $f$ is homogeneous, then $b_f(s) = b_{f,\bm 0}(s)$. This can be done using the natural $\mathbb C^*$-action.
		%
	\end{remark}
	\begin{remark}
		There are some generalizations of the Bernstein-Sato polynomial. It can be defined for an effective divisor on a smooth complex variety, see \cite{DMBG}. It can also be defined for a nonzero ideal $\mathfrak a \subseteq \mathbb C[\bm z]$, see \cite{Budur_Mustata_Saito_B_S_for_Ideal}. In addition, it can be generalized as an ideal in $\mathbb C[s_1,...,s_r]$ that defined for a tuple $(f_1,...,f_r)$ with $f_1,...,f_r \in \mathbb C[\bm z]\setminus \mathbb C$, see \cite{Zero_Loci_of_Bernstein-Sato_Ideals}. 
	\end{remark}

	\subsection{Archimedean zeta function}\label{Subsection_Archimedean_Zeta_Function}
	The definition and history of the archimedean zeta function was given in the introduction. Here, a more detailed review of this object is provided.

	Let $f\in \mathbb C[\bm z]\setminus \mathbb C$ and $b_f(s)$ be its $b$-function. Take $P\in A_n[s]$ such that $P\cdot f^{s+1} = b_f(s)\cdot f^s$. The proof of meromorphy in \cite{Bernstein-polynomial_Continuations} is sketched as follows.

	\vspace{0.5em}

	Suppose $P = \sum a_{\alpha,\bm \beta,\bm\gamma}s^\alpha \bm z^{\bm \beta} (\bm{\frac{\partial}{\partial z}})^{\bm \gamma}, a_{\alpha,\bm \beta,\bm \gamma} \in \mathbb C$. We define $P^* = \sum  a_{\alpha,\bm \beta,\bm\gamma}s^\alpha (-\bm{\frac{\partial}{\partial z}})^{\bm \gamma} \bm z^{\bm \beta}$. Further set $\overline P = \sum \overline{a_{\alpha,\bm \beta,\bm\gamma}} s^\alpha \bm{\bar z}^{\bm \beta} (\bm{\frac{\partial}{\partial \bar z}})^{\bm \gamma}$ and $\overline P^* = \sum \overline{a_{\alpha,\bm \beta,\bm\gamma}} s^\alpha (-\bm{\frac{\partial}{\partial \bar z}})^{\bm \gamma} \bm{\bar z}^{\bm \beta} $. Recall $\mathbb U = \{s \in \mathbb C \mid \mathrm{Re}\, s > 0\}$ in our notation.
	
	\vspace{0.3em}

	\noindent\textit{Step 1: $Z_{f,\varphi}(s)$ is holomorphic on $\mathbb U$.}
	
	This step can be done by the dominated convergence. The same spirit also appears in Subsection \ref{Subsection_Variation_of_Archimedean_Zeta_Function}. In fact, the argument tells us that integral expression is a holomorphic function of $s$ even on $\{s\in \mathbb C \mid \mathrm{Re}\, s > -\mathrm{lct}_{\{\varphi > 0\}}(f)\}$. Furthermore, the derivative is given by 
	\begin{displaymath}
		\frac{\mathrm{d}}{\mathrm{d}s} Z_{f,\varphi}(s) = \int_{\mathbb C^n} 2\vert f\vert^{2s} (\log \vert f\vert)\varphi\, \mathrm{d} \bm z.
	\end{displaymath}

	\noindent\textit{Step 2: Functional equation and extension.}
	
	Suppose $s\in\mathbb U$, then using integration by parts, we have
	\begin{align*}
		b_f(s) \overline{b_f(s)} \cdot Z_{f,\varphi}(s) & = \int_{\mathbb C^n}\big( Pf^{s+1} \cdot \overline P \bar f^{s+1}\big) \varphi\,  \mathrm{d} \bm z = \int_{\mathbb C^n} \vert f\vert^{2(s+1)} (P^*\overline P^*\varphi)\, \mathrm{d}\bm z.
	\end{align*}
	Notice that $P^*\overline P^*\varphi$ is a polynomial in $s$ with coefficients in $C_c^{\infty}(\mathbb C^n)$. Apply {Step 1}, then we see that the right-hand side is holomorphic on $\{s\in \mathbb C\mid \mathrm{Re}\, s > -1\}$. We thus conclude the proof by a simple induction.

	From the sketch, the following byproduct is obtained.
	\begin{lemma}\label{Fundamental_Lemma_of_Archi_B_Function} If $s_0$ is a pole of $Z_{f,\varphi}(s)$, then $s_0+\alpha$ is a root of $b_f(s)$ for some natural number $\alpha \leq \lfloor -\mathrm{Re}\, s_0\rfloor$. Moreover, if $\mathrm{Re}\, s_0 > -1$ and $s_0$ is an order $m$ pole, then $(s-s_0)^m \mid (b_f(s))^2$.
	\end{lemma}

	\begin{remark}\label{LCT_LARGEST_POLE}
		There is a simple observation that for $\varphi \geq 0$, $-\mathrm{lct}_{\{\varphi > 0\}}(f)$ is the largest pole of $Z_{f,\varphi}(s)$. Indeed, for real number $s > -\mathrm{lct}_{\{\varphi > 0\}}(f)$, $Z_{f,\varphi}(s)$ is a positive real number, and we have $ Z_{f,\varphi}(s) \to +\infty$ as $s \to (-\mathrm{lct}_{\{\varphi > 0\}}(f))^{+}$ by the defintion of log canonical thresholds.
	\end{remark}

	{Lemma \ref{Fundamental_Lemma_of_Archi_B_Function}} and {Theorem \ref{Ch1_Sec2_Kashiwara_Malgrange_Theorem}} are classical strategies in singularity theory. In \cite{Quelques}, Loeser proposed the converse problem of the first assertion in {Lemma \ref{Fundamental_Lemma_of_Archi_B_Function}}. That is, for $f\in\mathbb C[\bm z]\setminus \mathbb C$, does the following equality hold?
	\begin{displaymath}
		\bigcup_{\varphi \in C_c^{\infty}(\mathbb C^n)}\{\text{Poles of }Z_{f,\varphi}(s)\} = \{s_0-m\mid b_f(s_0) = 0, m \in \mathbb N\}.
	\end{displaymath}
	He verified the equality for plane curves and isolated quasi-homogeneous singularities. However, most recently, Davis, Lőrincz, and Yang provided a counterexample and thus answered the problem of Loeser in the negative, see \cite{DLY24}. A weaker version, in the sense of Kashiwara-Malgrange theorem, of this problem about monodromy eigenvalues was addressed positively. In \cite{Bar84}, Barlet proved that all monodromy eigenvalues are of the form $e^{2\pi \sqrt{-1} \alpha}$, where $\alpha$ is a pole of $Z_{f,\varphi}(s)$ for some $\varphi\in C_c^{\infty}(\mathbb C^n)$. In \cite{Bar86}, he further related the sizes of Jordan blocks of monodromies with the order of poles. In \cite{Lichtin89}, Lichtin refined the Kashiwara rationality theorem for the $b$-functions and gave a set of candidate roots.
	In \cite{ABCN17},  Artal Bartolo, Cassou-Nogu\`es, Luengo, and Melle-Hern\'andez  considered the residues and proved the Yano conjecture for the case in which the irreducible germ has two Puiseux pairs and its algebraic monodromy has distinct eigenvalues. The Yano conjecture was finally solved by Blanco, see \cite{Blanco31}.
	In \cite{Blanco19}, Blanco gave an explicit formula for residues of archimedean zeta functions for plane curves through resolutions.

	The archimedean zeta function can be generalized to a tuple of polynomials, see \cite{General_Archimeadean_Zeta}. There are also analogs of archimedean zeta functions, such as $p$-adic zeta functions and motivic zeta functions, see \cite{Igusa_Zeta_Function_Prologue} and \cite{Motivic_Igusa_Zeta_Function}. It is natural to ask parallel problems, such as Lemma \ref{Fundamental_Lemma_of_Archi_B_Function}, to these analogs. This provides a natural intuition for the monodromy conjecture.

	\subsection{Hyperplane arrangement and \texorpdfstring{$n/d$}{n/d}-conjecture}\label{Subsection_Hyperplane_Arrangement}
	
	In this subsection, we review some more conceptions related to hyperplane arrangements, especially results concerned with the $n/d$-conjecture.
	
	Let $A = \bigcup_{i=1}^r V_i \subseteq \mathbb C^{n}$ be a central hyperplane arrangement, i.e. each $V_i$ is a hyperplane in $\mathbb C^n$ with $\bm 0 \in V_i$. Suppose $V_i = \{f_i = 0\}$, where $f_1,...,f_r \in \mathbb C[\bm z] = \mathbb C[z_1,...,z_n]$ are homogeneous of degree one. Let $f = f_1^{b_1}\cdots f_r^{b_r},\bm b\in \mathbb N_{>0}^r$.
	 $A$ is called \textit{generic} if it is essential, indecomposable, and $\mathbb P(A)$ has simple normal crossings in $\mathbb P^{n-1}$, where $\mathbb P(A) = \{[\bm a] \in \mathbb P^{n-1} \mid \bm 0\neq \bm a \in A\}$. $f$ is called \textit{reduced} if $b_1 = ... =b_r = 1$.

	An \textit{edge} of $A$ is an intersection of some hyperplanes among $\{V_1,...,V_r\}$. Suppose $K$ is an edge of $A$ and $V_{l_1},...,V_{l_i}$ are all hyperplanes in $A$ containing $K$. Then $\overline A_K := \bigcup_{j = 1}^i (V_{l_j}/K) \subseteq \mathbb C^n/K$ defines a hyperplane arrangement $\mathbb C^n/K$. $K$ is called \textit{dense} if $\overline A_K$ is indecomposable.

	
	\vspace{0.5em}

	Next, we review the \textit{canonical log resolution} of $(\mathbb C^n ,A)$. We refer to \cite{Cohomology_of_Local_System_on_the_Complement_of_Hyperplane_Arrangement} and \cite{Jumping_Coefficient_and_Spectrum_of_a_Hyperplane_Arrangement_Budur_Saito} for more details.

	Let $\mathcal S$ be the set of edges. Then the canonical log resolution of $(\mathbb C^n,A)$ is given by blowing up strict transforms of edges with dimensions increasing. The canonical log resolution is denoted by $\mu : Y \to \mathbb C^{n}$. The strict transforms and exceptional divisors are indexed by $\mathcal S$. So, in the resolution $\mu$, for an edge $W\in \mathcal S$, we denote the divisor from $W$ by $E_W$. 
	
	There is another canonical log resolution given by blowing up strict transforms of all dense edges with dimensions increasing, see \cite{Monodromy_Conjecture_for_Hyperplane_Arrangement_Budur_Mustata_Teitler}.
	
	\vspace{0.5em}

	
	In the following part, we list some existing results related to the $n/d$-conjecture.
	
	The case when $f$ is reduced and $\{f =0\}$ is generic can be covered by the formula for the whole $b_f(s)$ given in \cite{B_S_for_Centeral_Generic_Hyperplane_Arrangements_Walther} and  \cite{B_S_for_Centeral_Generic_Hyperplane_Arrangements_Saito}.
	\begin{theorem}\label{Bernstein_for_Reduced_generic}
		If $f$ is reduced and $A = \{f = 0\}$ is generic, then the Bernstein-Sato polynomial of $f$ is
		\begin{displaymath}
			b_f(s) = (s+1)^{w} \prod_{j = 0}^{2d-n-2} (s+\frac{j+n}{d}),
		\end{displaymath}
		where Walther gave the formula and proved $w \in \{n-1,n-2\}$ in \cite{B_S_for_Centeral_Generic_Hyperplane_Arrangements_Walther}  and Saito further proved $w = n-1$ in \cite{B_S_for_Centeral_Generic_Hyperplane_Arrangements_Saito}. 
	\end{theorem}
	
	In \cite{Local_Zeta_Function_And_b_Function_of_Certain_Hyperplane_Arragement_Budur_Saito_Sergey}, Budur, Saito, and Yuzvinsky specifically studied the $n/d$-conjecture. Some interesting results were provided.
	\begin{theorem} $n/d$-conjecture holds in the following cases.
		
		\noindent(1) $\bm 0$ is a good dense edge of $f$.
		
		\noindent(2) $f$ is reduced, and $n \leq 3$.
		
		\noindent(3) $f$ is reduced, $(n,d) = 1$, and $V_r$ is generic relative to $V_1,...,V_{r-1}$. 
	\end{theorem}
	\begin{remark}
		Case (1) is equivalent to saying that ${n}/{d}$ is the log canonical threshold of $f$. $V_r$ is generic relative to $V_1,...,V_{r-1}$ means that any nonzero interaction of $V_j$ with $j< r$ is not contained in $V_r$.
	\end{remark}
	
	In \cite{Tame_Arrangement}, Bath computed the roots of the $b$-function for tame hyperplane arrangements. \cite[Theorem 1.3]{Tame_Arrangement} covers the $n/d$-conjecture for tame arrangements.
	\begin{corollary}
		Suppose $A = \{f = 0\}$ is a central essential indecomposable hyperplane arrangement and $f$ is tame, then the $n/d$-conjecture holds for $f$.
	\end{corollary}
	\begin{remark}\label{tame leq 3}
		Here $f$ is called \textit{tame} if  for each $k\in \mathbb N_{>0}$ and $x\in \mathbb C^n$, the $k$-th logarithmic differential $\Omega_{\mathbb C^n,x}^k(\log(\mathrm{div}\, f))$ has projective dimension $\leq k$. By  \cite[{Remark 2.5}(3)]{Tame_Arrangement}, if $n \leq 3$, then $f$ is always tame. Hence, the $n/d$-conjecture holds when $n\leq 3$.
	\end{remark}
	
	In \cite{MPVIHA}, Budur and the authors studied the motivic principal value integral for canonical resolutions of hyperplane arrangements and proved a partial converse of the second assertion of {Theorem \ref{Weak Monodromy Conjecture for hyperplane arrangements}}. It indicated that the strong monodromy conjecture implies the $n/d$-conjecture for generic multiplicities as follows.
	\begin{corollary}\label{propC2} 
		Consider a central essential indecomposable arrangement in $\mathbb C^n$ of $d$ mutually distinct hyperplanes defined by linear polynomials $f_i\in\mathbb C[x_1,...,x_n]$ with $i=1,..., d$. 
		There exists a nonzero product $\Theta$ of polynomials of degree one in $d$ variables such that if $\bm m \in \mathbb Z_{>0}^d\setminus \{\Theta = 0\}$ then $-n/\sum_{i=1}^dm_i$ is a pole of the motivic zeta function of $g=f_1^{m_1}\cdots f_d^{m_d}$. 
		In particular, for these $g$, the strong monodromy conjecture implies the $n/d$-conjecture. 
	\end{corollary}
	\begin{remark}
		However, $-n/d$ is not always a pole of $Z_{g}^{\mathrm{mot}}(s)$, cf. the example given by Veys in \cite[Appendix]{Local_Zeta_Function_And_b_Function_of_Certain_Hyperplane_Arragement_Budur_Saito_Sergey}.
	\end{remark}

	\section{Variation of archimedean zeta function}\label{section3}

	Before starting this section, we fix some notation.
	
	Fix $f_1,...,f_r \in \mathbb C[\bm z]\setminus \mathbb C = \mathbb C[z_1,...,z_n]\setminus \mathbb C$.  Let $\mu: Y \to \mathbb C^n$ be a log resolution of $f = f_1f_2\cdots f_r$. Suppose $\mu^*( \mathrm{div}\, f) = \sum_{j\in I} N_j E_j$ and $K_{\mu} = \sum_{j\in I} (\nu_j-1) E_j$, where $N_{j} > 0$ for all $j \in I$. Here, $\{E_j\}_{j\in I}$ is a collection of simple normal crossings divisors on $Y$. Suppose $\mu^*(\mathrm{div} f_i) = \sum_{j\in I} N_{ij} E_j$. For $j\in I$, we set $P_j : \mathbb C^r \to \mathbb C,\ \bm b =(b_1,...,b_r) \mapsto \sum_{i=1}^r b_iN_{ij}$. For $I' \subseteq I$, denote by $E_{I'}$ the intersection $\bigcap_{i\in I'} E_i$.
	
	Since $\mathrm{Supp}\, \varphi$ is compact and $\mu$ is proper, $\mu^{-1}(\mathrm{Supp}\, \varphi)$ is also compact. Then we can cover $\mu^{-1}(\mathrm{Supp}\, \varphi)$ by a finite collection of open subsets of $Y$, say $\{U_k\}_{k\in J}$, such that (a) for each $k$, there are coordinates $x_{k1},...,x_{kn}$ satisfying that for each $j$, $E_j \cap U_k = \{x_{kl} = 0\}$ for some $l$ or empty; (b) the coordinates $x_{k1},...,x_{kn}$ map $U_k$ biholomorphically to a polydisk in $\mathbb C^n$ centering at the origin.
	
	For each $U_k$, we assume that $E_{j_p}\cap U_k = \{x_{kp} = 0\}$ for $p = 1,...,v_k$ and $E_j\cap U_k = \varnothing$ for other divisors. Let $\{\phi_k\}_{k\in J}$ be a partition of unity with respect to $\{U_k\}_{k\in J}$.

	\subsection{Definition and meromorphy}	\label{Subsection_Variation_of_Archimedean_Zeta_Function}
	
	In this subsection, we will define the \textit{variation of archimedean zeta function}, prove its meromorphy, and provide its candidate poles. We start with some definition and a technical lemma.

	\begin{definition}
		Let $Z$ be a Lebesgue null set in $\mathbb C^n$. A Lebesgue-measurable function $\psi : (\mathbb C^{n}\setminus Z) \times \mathbb C^l\to \mathbb C,\ (\bm z,\bm x) \mapsto \psi(\bm z,\bm x)$ is called \textit{suitable} if the following hold.
		
		\noindent(S1) $\mathrm{Supp}\, \psi \subseteq K \times \mathbb C^l$ for some compact subset $K\subseteq \mathbb C^n$.
		
		\noindent(S2) For an arbitrary bounded subset $K' \subseteq \mathbb C^l$, $\psi : (\mathbb C^n\setminus Z) \times K' \to \mathbb C$ is bounded.  
		
		\noindent(S3) For all $\bm z_0 \in \mathbb C^n\setminus Z'$, $\psi(\bm z_0,-)$ is holomorphic. Moreover, for all $k=1,...,l$, $\frac{\partial \psi}{\partial x_k} : (\mathbb C^n\setminus Z) \times \mathbb C^l \to \mathbb C$ satisfies (S2). 
	\end{definition}
	
	\begin{example}\label{key examples}
		(1) Let $U\subseteq \mathbb C^n$ be an open subset and $g_1,...,g_{l-1} : U \to \mathbb C$ are holomorphic functions. $\varphi$ is a bounded measurable function on $\mathbb C^n$, supported in $U$. Assume all $g_1,...,g_{l-1}$ are nowhere zero on $U$, then $\vert g_1\vert^{2x_1x_{l}}\cdots\vert g_{l-1}\vert^{2x_{l-1}x_{l}} \varphi : \mathbb C^{n} \times \mathbb C^{l} \to \mathbb C$ is suitable.
		
		\noindent(2) If $\varphi \in C^{\infty}_c(\mathbb C^n) = C_c^{\infty}(\mathbb C^n,\mathbb C)$, then $\Phi : (\bm z,\bm x) \to \varphi(\bm z)$ is suitable. Furthermore, let $z_i = \rho_i e^{\theta_i\sqrt{-1}}$ be the presentation by polar coordinates, then for all $\bm m \in \mathbb N^n$, $\frac{\bm \partial^{\vert \bm m\vert} \Phi}{\bm{\partial \rho}^{\bm m}} : (\mathbb C^n\setminus \{\bm 0\}) \times \mathbb C^l \to \mathbb C$ is also suitable.
		
		\noindent(3) The product of two suitable functions defined on $(\mathbb C^n\setminus Z) \times \mathbb C^l$ is again suitable.
	\end{example}

	Take a suitable function $\psi : (\mathbb C^n\setminus Z) \times \mathbb C^{l} \to \mathbb C$. We introduce the following auxiliary function $L_{\psi}$. Similar to the classical archimedean zeta function, it is also expressed as an integral. Recall that $\mathrm{d} \bm z = (\frac{\sqrt{-1}}{2})^n \mathrm{d}z_1\wedge \mathrm{d} \bar z_1\wedge \cdots \wedge \mathrm{d}z_{n} \wedge \mathrm{d} \bar z_{n}$ is the normalized measure.
	\begin{displaymath}
		L_{\psi}(\bm a,\bm x) := \int_{\mathbb C^n} \prod_{i=1}^r \vert f_i(\bm z)\vert^{2a_i}\psi(\bm z,\bm x)\, \mathrm{d} \bm z,\ (\bm a,\bm x)\in \mathbb U^r\times \mathbb C^l. 
	\end{displaymath}
	\begin{lemma}\label{main technique}
		$L_{\psi}(\bm a,\bm x)$ is holomorphic on $\mathbb U^r \times \mathbb C^l$.
	\end{lemma}	
	\begin{proof}
		The problem is local. So we fix a tuple $(\bm a,\bm x) \in \mathbb U^r \times \mathbb C^{l}$ and consider $L_\psi$ in its neighborhood. Let $K\subseteq \mathbb C^n$ be a compact subset such that $\mathrm{Supp}\, \psi \subseteq K \times \mathbb C^{l}$. Since the result does not change if we scale $f_1,...,f_r$ by a nonzero complex number, we may assume $\vert f_i(\bm z)\vert < 1$ for all $\bm z\in K$.

		Since $\mathrm{Re}\, a_i > 0$ for all $i$, there exists $\alpha > 0$ such that $\mathrm{Re}\, a_i > \alpha$ for all $i$. Then we have $\prod_{i=1}^r\vert f_i(\bm z)\vert^{2a_i}\vert \psi(\bm z,\bm x)\vert < \vert \psi(\bm z,\bm x)\vert $ is integrable on $\mathbb C^n$. Hence, $L_{\psi}$ is well-defined. We will show its continuity and then its holomorphy. 
		
		Take $\varepsilon \in (0,\alpha/2)$ such that $\mathbb B((\bm a,\bm x) , 2\varepsilon) \subset \mathbb U^{r} \times \mathbb C^{l}$, where $\mathbb B(p,\varepsilon)$ denotes the Euclidean ball. 
		
		\noindent\textit{Step 1: Continuity.}
		
		For $(\bm b,\bm y) \in \mathbb C^{r+l}$, with $\vert (\bm b,\bm y)\vert < \varepsilon$ we have
		\begin{align*}
			&\vert L_{\psi}(\bm a+\bm b,\bm x+\bm y)-L_{\psi}(\bm a,\bm x)\vert  \leq  \int_{\mathbb C^n} \bigg( \vert f_i\vert^{2\mathrm{Re}(a_i+b_i)} \vert \psi(\bm z,\bm x+\bm y) - \psi(\bm z,\bm x)\vert \\
			&+ \sum_{j=1}^r\prod_{i=1}^{j-1} \vert f_i\vert^{2\mathrm{Re}\, a_i}\cdot \vert f_j\vert^{2\mathrm{Re}\, a_j}\cdot \big\vert \vert f_j\vert^{2b_j}-1\big\vert\cdot \prod_{i=j+1}^r\vert f_i\vert^{2\mathrm{Re}(a_i+b_i)} \cdot \vert\psi(\bm z,\bm x)\vert \bigg) \, \mathrm{d} \bm z .
		\end{align*}
		Since $\mathrm{Re}(b_i) > -\alpha$, we have 
		\begin{displaymath}
			\vert f_i\vert^{2\mathrm{Re}\, a_i}\cdot \big\vert \vert f_i\vert^{2b_i}-1\big\vert =\big\vert \vert f_i\vert^{2\mathrm{Re}a_i} - \vert f\vert^{2\mathrm{Re}\, a_i+2b_i} \big\vert \leq \vert f_i\vert^{2\mathrm{Re}\, a_i} + \vert f\vert^{2\mathrm{Re}\, a_i+2\mathrm{Re}\, b_i} < 1+1 = 2
		\end{displaymath}
		for all $i$. Moreover, $\big\vert \vert f_i\vert^{2\mathrm{Re}a_i} - \vert f\vert^{2\mathrm{Re}\, a_i+2b_i} \big\vert$ tends to zero as $\bm b \to 0$. Hence, take $M = \sup_{\bm x' \in \mathbb B(\bm x,\varepsilon)} \Vert \psi(-,\bm x')\Vert_{L^{\infty}}$, then $(2n+2)M \cdot \bm{1}_K$ is a dominating function, and continuity is verified by using the dominated convergence. 
		
		\noindent\textit{Step 2: Holomorphy.}
		
		Consider $\frac{\partial L}{\partial a_1}$ first. For a complex number $b_1$ with $0< \vert b_1\vert < \varepsilon$, we have
		\begin{displaymath}
			\frac{L_{\psi}(a_1+b_1,\bm a',s)-L_{\psi}(a_1,\bm a',s)}{b_1} = \int_{\mathbb C^n} \prod_{i=1}^r \vert f_i\vert^{2a_i} \frac{\vert f_1\vert^{2b_1}-1}{b_1} \psi\, \mathrm{d} \bm z.
		\end{displaymath}
		where $\bm a' = (a_2,...,a_r)$. For $\bm z \not\in \{f_1 = 0\}$, by the Lagrange inequality, we have
		\begin{displaymath}
			\vert \frac{\vert f_1(\bm z)\vert^{2b_1}-1}{b_1}\vert \leq \big\vert \log \vert f_1(\bm z)\vert \big\vert\cdot \vert f_i(\bm z)\vert^{2\mathrm{Re}\,\xi} \leq  \big\vert \log \vert f_1(\bm z)\vert \big\vert\cdot \vert f_i(\bm z)\vert^{-2\varepsilon},
		\end{displaymath}
		where $\xi = \xi_0 b_1$ for some $\xi_0 \in [0,1]$. The last inequality holds since $\vert f(\bm z)\vert < 1$ and $\mathrm{Re}\, \xi \geq -\vert b_1\vert \geq -\varepsilon$. Notice that for $\beta >0$, $u\mapsto u^{\beta} \log u$ is a bounded function on $(0,1)$. Since $\{f_1 = 0\}$ is a null subset of $\mathbb C^n$ and $\mathrm{Re}\, a_i - \varepsilon > \alpha-\varepsilon > 0$, $\prod_{i=2}^r\vert f_i\vert^{2\mathrm{Re}\, a_i} \cdot \vert f_1\vert^{2\mathrm{Re}(a_1)-2\varepsilon} \big\vert \log \vert f_1\vert \big\vert \cdot \vert\psi\vert$ is a dominating function. Therefore, again by the dominated convergence, we have
		\begin{displaymath}
			\frac{\partial L_{\psi}}{\partial a_1}(\bm a,\bm x) = \lim_{\vert b_1\vert < \varepsilon, b_1 \to 0} \frac{L_{\psi}(a_1+b_1,\bm a',\bm x)-L_{\psi}(a_1,\bm a',\bm x)}{b_1} = 2\int_{\mathbb C^n} \prod_{i=1}^r\big( \vert f_i\vert^{2a_i} \log \vert f_1\vert \big) \psi\, \mathrm{d} \bm z.
		\end{displaymath}
		Similarly, all $\frac{\partial L_{\psi}}{\partial a_2},...,\frac{\partial L_{\psi}}{\partial a_r}$ exist. The existence of $\frac{\partial L_{\psi}}{\partial x_j}$ follows from (S3) and a similar argument using the dominated convergence.
	\end{proof}
	
	\vspace{0.5em}
	
	We apply some changes of variables to $L_{\psi}$. Let $\lambda \in \mathbb R_{>0}$ and we associate $\lambda$ with a half-strip as follows.
	\begin{displaymath}
		\mathbb H_{\lambda} := \{z \in \mathbb C \mid \mathrm{Re}\, z > \lambda, \vert \mathrm{Im}\, z\vert < \lambda\}.
	\end{displaymath}	
	We have a natural map:
	\begin{displaymath}
		\mathbb H_{\lambda}^r \times  \mathbb H_{\lambda} \to \mathbb U^r, (\bm a,s) \mapsto (a_1s,...,a_rs).
	\end{displaymath}
	Let $\varphi \in C_c^{\infty}(\mathbb C^n)$, then by {Lemma \ref{main technique}}, the following map is holomorphic.
	\begin{displaymath}
		Z_{\varphi} : \mathbb H_{\lambda}^r\times  \mathbb H_{\lambda} \to \mathbb C,\ (\bm a,s) \mapsto \int_{\mathbb C^n} \prod_{i=1}^r \vert f_i(\bm z)\vert^{2a_is}\varphi(\bm z)\, \mathrm{d} \bm z.
	\end{displaymath}
	
	Quite similar to classical archimedean zeta functions, we can extend $Z_{\varphi}$ meromorphically to $\mathbb C^{r+1}$. The strategy also involves a log resolution.
	
	\begin{definition}
		
		The notation is as at the beginning of Section \ref{section3}. Let $F: \mathbb C^{r+1}\to \mathbb C$ be a holomorphic function. We call $I' \subseteq I$ \textit{$F$-possible}, if $E_{I'} \neq \varnothing$ and there are natural numbers $\{\beta_j\}_{j\in I'}$ such that functions $P_j(-)+\nu_j+\beta_j/2,j\in I'$ are all complex multiples of $F$. 
		Let $m(F)$ be the largest cardinality of an $F$-possible subset of $I$.
	\end{definition}
	
	The following proposition denotes the meromorphic nature of $Z_{\varphi}$.
	\begin{proposition}\label{meromorphic extensino of VAZF}
		The notation is as at the beginning of Section \ref{section3}. $Z_{\varphi}$ admits a meromorphic extension to $\mathbb C^{r+1}$. Furthermore, we have the following.
		
		\noindent(1) The polar set of $Z_{\varphi}$ is contained in 
		\begin{displaymath}
			\Lambda := \bigcup_{\beta \in \mathbb N}\bigcup_{j\in I} \{(\bm a,s) \in \mathbb C^{r+1} \mid P_j(\bm a)s+\nu_j+\beta/2 = 0\} \subseteq \mathbb C^{r+1}.
		\end{displaymath}
		
		\noindent(2) For a nonconstant holomorphic function $F : \mathbb C^{r+1}\to \mathbb C, (\bm a,s) \mapsto \sum_{i=1}^r c_ia_is + c_0$ with $c_0,...,c_r \in \mathbb C$, $\{F = 0\}$ is not contained in the polar set of $F^{m(F)} \cdot Z_{\varphi}$. That is, the order of $\{F = 0\}$ as a polar component not greater than $m(F)$.
	\end{proposition}
	\begin{proof}
		The notation is as at the beginning of Section \ref{section3}. By the change of variable formula, we have
		\begin{displaymath}
			Z_{\varphi}(\bm a,s) = \sum_{k\in J} \int_{\mathbb C^n} \prod_{p=1}^{v_k} \vert x_{kp}\vert^{2(P_{j_p}(\bm a)s+\nu_{j_p}-1)} \big(\prod_{i=1}^r \vert u_{ki}\vert^{2a_is}  \big) \vert u_k\vert^2(\varphi\circ \mu) \phi_k\, \mathrm{d} \bm x_k, \tag{\ref{meromorphic extensino of VAZF}.1}
		\end{displaymath}
		where we set $(f_i\circ \mu)|_{U_k} = \prod_{p=1}^{v_k}x_{kp}^{N_{ij_p}} \cdot u_{ki}$ and $\det \big(\mathrm{Jac}_\mu\big) |_{U_k} = \prod_{p=1}^{v_k} x_{kp}^{\nu_{j_p}-1}\cdot u_k$. Hence, $u_{k},u_{k1},...,u_{kr}$ are all holomorphic and everywhere nonzero on $U_k$ (identified with a polydisk in $\mathbb C^n$). Let $Z_k(\bm a,s)$ denote the integral in $k$-place of the right-hand side of (\ref{meromorphic extensino of VAZF}.1). It suffices to prove the meromorphy, (1), and (2) for each $Z_k$. The argument is similar to the one for the classical archimedean zeta function.
		
		We drop $k$ in the notation and the target object is
		\begin{displaymath}
			Z(\bm a,s) = \int_{\mathbb C^n} \prod_{p=1}^{v} \vert x_{p}\vert^{2(P_{j_p}(\bm a)s+\nu_{j_p}-1)} \big(\prod_{i=1}^r \vert u_{i}\vert^{2a_is}  \vert u\vert^2(\varphi\circ \mu) \phi\big)\, \mathrm{d} \bm x.
		\end{displaymath}
		Consider $(\bm a,s) \in \mathbb H_\lambda^{r+1}$. Let $z_i = \rho_i e^{\theta_i\sqrt{-1}}$ be the polar coordinate, then
		\begin{displaymath}
			Z(\bm a,s) = \int_{\mathbb R_+^n}\int_{[0,2\pi]^n} \prod_{p=1}^{v}\rho_p^{2(P_{j_p}(\bm a)s+\nu_{j_p})-1}  \big(\prod_{i=1}^r \vert u_{i}\vert^{2a_is}  \vert u\vert^2(\varphi\circ \mu) \phi\big) \, \mathrm{d}\bm \rho \mathrm{d} \bm \theta.
		\end{displaymath}
		Integrating by parts to $\beta_1,...,\beta_v$, we have
		\begin{displaymath}
			Z(\bm a,s) = (-1)^{\vert\bm  \beta\vert}C_{\bm \beta}\int_{\mathbb R_+^n}\int_{[0,2\pi]^n} \prod_{p=1}^v\rho_p^{2(P_{j_p}(\bm a)s+\nu_{j_p})+\beta_p-1}  \frac{\bm{\partial^{\vert  \beta\vert}}}{\bm{\partial \rho^{ \beta}}}\big(\prod_{i=1}^r \vert u_{i}\vert^{2a_is}  \vert u\vert^2(\varphi\circ \mu) \phi\big) \, \mathrm{d}\bm \rho \mathrm{d} \bm \theta, \tag{\ref{meromorphic extensino of VAZF}.2}
		\end{displaymath}
		for all $\bm \beta \in \mathbb N^v\times \{0\}^{n-v}$. Here $\displaystyle C_{\bm \beta} = \prod_{p=1}^v\prod_{t=0}^{\beta_i-1} \frac{1}{2(P_{j_p}(\bm a)s+v_{j_p})+t}$. 
		
		By {Example \ref{key examples}} and {Lemma \ref{main technique}}, the right-hand side of (\ref{meromorphic extensino of VAZF}.2) is holomorphic in the region
		\begin{displaymath}
			D_{\bm \beta} := \bigcap_{p=1}^v \{(\bm a,s)\in \mathbb C^{r+1} \mid \mathrm{Re}\big(2(P_{j_p}(\bm a)s+\nu_{j_p}) +\beta_p-2\big) > 0\}.
		\end{displaymath}
		Since $\bigcup_{\bm \beta \in \mathbb N^v\times \{0\}^{n-v}} D_{\bm \beta} = \mathbb C^{r+1}$, we conclude the proof of meromorphy. (1) and (2) follow straightforwardly from the process.
	\end{proof}
	
	\begin{definition}
		For $f_1,...,f_r \in \mathbb C[\bm z]\setminus \mathbb C$ and $\varphi \in C_c^{\infty}(\mathbb C^n)$, we define the \textit{variation of archimedean zeta function} $Z_{\varphi}(\bm a,s)$ to be the meromorphic extension of 
		\begin{displaymath}
			\mathbb H_{\lambda}^r\times  \mathbb H_{\lambda} \to \mathbb C,\ (\bm a,s) \mapsto \int_{\mathbb C^n} \prod_{i=1}^r \vert f_i(\bm z)\vert^{2a_is}\varphi(\bm z)\, \mathrm{d} \bm z
		\end{displaymath}
		to the whole $\mathbb C^{r+1}$. It is clear that the definition of $Z_{\varphi}$ is independent of $\lambda$.
	\end{definition}
	\begin{remark}
		The variation of the archimedean zeta function will be a useful approach to studying the poles of the classical archimedean zeta function. For a fixed $\bm a_0 \in \mathbb C^r$, $\Lambda \cap \big(\{\bm a_0\} \times \mathbb C\big)$ is a discrete set $C_{\bm a_0} \in \mathbb C$. In particular, if $\bm a_0 \in \mathbb N_{>0}^r$, then $Z_{\varphi}(\bm a_0,-) = Z_{\bm f^{\bm a_0},\varphi}(-)$ on $\mathbb C\setminus C_{\bm a_0}$. This is because the two holomorphic functions coincide on $\mathbb U = \{s \in \mathbb C\mid \mathrm{Re}\, s > 0\}$. So, the variation version can specialize to the classical one with no information lost.
		
		However, one should be careful when trying to restrict a meromorphic function to a submanifold. To circumvent the subtleties, our tactic is to multiply $Z_{\varphi}$ by some holomorphic functions and make the new function holomorphic on some neighborhoods of some points in the polar set of $Z_{\varphi}$. More precisely, if $(\bm a_0,s_0) \in \mathbb N_{>0}^r \times \mathbb C$ is in the polar set of $Z_{\varphi}$, let $F : \mathbb C^{r+1}\to \mathbb C$ be a holomorphic function such that $F\cdot Z_{\varphi}$ is holomorphic around $(\bm a_0,s_0)$. Then we can restrict $F\cdot Z_{\varphi}$ to a neighborhood of $\{\bm a_0\} \times \{s_0\}$ in $\{\bm a_0\} \times \mathbb C$, and we have $(F\cdot Z_{\varphi})(\bm a_0,s_0) = \big(F(\bm a_0,-)\cdot Z_{\bm f^{\bm a_0},\varphi}\big)(s_0)$.
		%
	\end{remark}
	
	Moreover, $Z_{\varphi}$ admits a symmetry under scaling.
	\begin{lemma}\label{symmerty under scaling}
		For all $b \in \mathbb C^*$, $Z_\varphi({b\cdot \bm a, s}) = Z_\varphi(\bm a,b\cdot s)$ as meromorphic functions of $(\bm a,s)$. We call this symmetry \textit{self-scaling}. 
	\end{lemma}
	\begin{proof}
		From Proposition \ref{meromorphic extensino of VAZF}, we deduce that the polar sets of $Z_{\varphi}(b\cdot -,-)$ and $Z_{\varphi}(-,b\cdot -)$ are both contained in
		\begin{displaymath}
			\Lambda_{b} := \bigcup_{\beta \in \mathbb N}\bigcup_{j\in I} \{(\bm a,s) \in \mathbb C^{r+1} \mid b\cdot P_j(\bm a)s+\nu_j+\beta/2 = 0\} \subseteq \mathbb C^{r+1},
		\end{displaymath}
		which is nowhere dense in $\mathbb C^{r+1}$. So, it suffices to prove the equality each fixed $(\bm a,s) \in \mathbb C^{r+1} \setminus \Lambda_{b}$. 
		
		We consider the proof of Proposition \ref{meromorphic extensino of VAZF} and adopt the notation there. Since $Z_{\varphi}(\bm a,s) = \sum_{k\in J} Z_k(\bm a,s)$, we only need to prove the equality for each $Z_k$. Again, drop $k$ in the notation. Take $\bm \beta \in \mathbb N_{>0}^v \times \{0\}^{n-v}$ such that $(b\cdot \bm a,s) \in D_{\bm \beta}$, then $(\bm a,b\cdot s) \in D_{\bm \beta}$ according to the definition of $D_{\bm\beta}$. In the sequel, $Z(b\cdot \bm a,s) = Z(\bm a,b\cdot s)$ follows immediately from (\ref{meromorphic extensino of VAZF}.2).
		%
	\end{proof}

	\subsection{Generic survival of poles and $n/d$-conjecture for generic multiplicities}
	
	We still adopt the notation at the beginning of Section \ref{section3} and Notation \ref{N&C}.

	For a non-constant holomorphic function $F: (\bm a,s) \mapsto \sum_{i=1}^r c_i a_is + c_0,\ c_0,...,c_r\in \mathbb C,c_0\neq 0$, let $g_{_F}: \bm a\mapsto  \sum_{i=1}^r c_ia_i$. We identify $V_F := \mathbb C^r\setminus \{g_{_F} = 0\}$ with $\{F=0\}$ through the isomorphism $\theta : V_F \overset{\simeq}{\longrightarrow} \{F = 0\} \subseteq \mathbb C^{r+1}, \bm a \mapsto (\bm a,-\frac{c_0}{g_{_F}(\bm a)})$ 
	
	\begin{definition}
		For a non-constant holomorphic function $F: \mathbb C^{r+1} \to \mathbb C$, we say that $(Z_{\varphi},\mu)$ is \textit{of order one} at $F$ if $m(F) = 1$ and the zero set of $\widetilde Z_{\varphi,F} := F\cdot Z_{\varphi}$ does not contain $\{F = 0\}$. 
	\end{definition}
	\begin{remark}\label{XFphi}
		From the definition, the preimage of the union of zero and polar sets of $\widetilde Z_{\varphi,F}$ through $\theta$, denoted by $X_{F,\varphi}$, is analytically closed in $V_F$ and nowhere dense.
	\end{remark}

	\begin{lemma}\label{variation pole lemma}
		Suppose $(Z_\varphi,\mu)$ is of order one at $F = \sum_{i=1}^r c_ia_is+c_0,\ c_0,...,c_r\in \mathbb N_{\geq 0}$ for all $j$.  Then $\widetilde Z_{\varphi,F}$ is also self-scaling. Consequently, $X_{F,\varphi}$ is a cone without the vertex in $\mathbb C^r$ and $
		X_{F,\varphi,1} := X_{F,\varphi} \cap \{g_{_F} = 1\} \subseteq \{g_{_F} = 1\} \simeq \mathbb C^{r-1}$ is analytically closed and nowhere dense. 
	\end{lemma}
	\begin{proof}
		The self-scaling of $\widetilde Z_{\varphi,F}$ follows immediately from the self-scaling of $F$ and $Z_\varphi$, which is shown in {Lemma \ref{symmerty under scaling}}. The rest follows from {Remark \ref{XFphi}}. 
	\end{proof}
	
	Now we specialize the variation of archimedean zeta function to the classical one. 
	

	\begin{definition}
		We call $\bm a \in \mathbb N_{>0}^r$ a \textit{good tuple} at $j\in I$ if $-\frac{P_{j'}(\bm a)}{P_{j}(\bm a)}\nu_j+\nu_{j'} > 0$ for all $j'\neq j$. In particular, $\frac{\nu_j}{P_j(\bm a)}$ is the log canonical threshold of $f$.
	\end{definition}
	\begin{proof}[Proof of Theorem \ref{generic existence}]
		Take $F : (\bm a,s)\mapsto P_j(\bm a)s+\nu_j$. By {Lemma \ref{variation pole lemma}}, it suffices to prove $(Z_\varphi,\mu)$ is of order one at $F$. Since $P_{j'}(\bm a_0)\big(-\frac{\nu_j}{P_{j}(\bm a_0)}\big) + \nu_{j'} + \beta/2 > 0$ for all $j' \neq j$ and $\beta \geq 0$, $P_{j'}(-)s+\nu_{j'}+\beta/2$ cannot be a complex multiple of $F$. Hence, we have $m(F) = 1$.
		
		Let $\widetilde Z_{\varphi,F}(\bm a,s) := F\cdot Z_{\varphi}(\bm a,s)$ be as before. Then $\{F = 0\}$ is not in the polar set of $\widetilde Z_{\varphi,F}$. Since $(\bm a_0,-\frac{\nu_j}{P_j(\bm a_0)})$ does not belong to any other components of the polar set of $Z_{\varphi}$, $\widetilde Z_{\varphi,F}$ is holomorphic in a neighborhood of $(\bm a_0,-\frac{\nu_j}{P_j(\bm a_0)})$. It suffices to show $\widetilde Z_{\varphi,F}(\bm a_0,-\frac{\nu_j}{P_j(\bm a_0)}) \neq 0$. Let $\widetilde Z_{\bm f^{\bm a_0},\varphi,F}(s) = (s+\frac{\nu_j}{P_j(\bm a_0)})Z_{\bm f^{\bm a_0},\varphi}(s)$. Since $\widetilde Z_{\varphi,F}|_{\{\bm a_0\} \times (\mathbb C\setminus C_{\bm a_0})}$ and $\widetilde Z_{\bm f^{\bm a_0},\varphi,F}$ are holomorphic and coincide on $\mathbb C \setminus C_{\bm a_0}$, we deduce that $\widetilde Z_{\bm f^{\bm a_0},\varphi,F}(-\frac{\nu_j}{P_j(\bm a_0)}) = \widetilde Z_{\varphi,F}(\bm a_0,-\frac{\nu_j}{P_j(\bm a_0)})$ is finite. So $-\frac{\nu_j}{P_j(\bm a_0)}$ is a candidate pole of $Z_{\bm f^{\bm a_0},\varphi}$ of order $\leq 1$.
		
		Moreover, according to our assumption, $\frac{\nu_j}{P_j(\bm a_0)}$ is the log canonical threshold of $\bm f^{\bm a_0}$. In addition, since $\varphi \geq 0$ and $\varphi|_{\mu(E_j)} \neq 0$, $-\frac{\nu_j}{P_j(\bm a_0)}$ is the largest pole of $Z_{\bm f^{\bm a_0},\varphi}(s)$. Therefore, we have $\widetilde Z_{\varphi,F}(\bm a_0,-\frac{\nu_j}{P_j(\bm a_0)}) = \widetilde Z_{\bm f^{\bm a_0},\varphi,F}(-\frac{\nu_j}{P_j(\bm a_0)}) \neq 0$.
	\end{proof}

	Then we apply our result to hyperplane arrangements. 
	
	Suppose $f_1,...,f_r \in \mathbb C[\bm z]$ are homogeneous polynomials of degree one and are pairwise linearly independent. Let $V_i := \{f_i = 0\}$, and we further assume that $A := \{f_1f_2\cdots f_r = 0\}$ is a central essential indecomposable hyperplane arrangement.
	
	Let $\mu : Y \to \mathbb C^n$ be the canonical log resolution of $f$, achieved by blowing up the strict transforms of all edges with dimensions increasing. Let $\mathcal S$ be the set of edges and $E_W$ be the corresponding strict transform or exceptional divisor of $W \in \mathcal S$. We therefore replace the notation $P_j,N_{ij},\nu_j,j\in I$ by $P_W,N_{iW},\nu_{W},W\in \mathcal S$. We shall construct a good tuple $\bm c$ at $\bm 0\in \mathcal S$.
	
	The numerical information of $\mu$ is given by $N_{lW} = \begin{cases} 1, & W \subseteq V_l\\ 0, & \text{else}\end{cases}$ and $\nu_W = \mathrm{codim} W$, see \cite{Cohomology_of_Local_System_on_the_Complement_of_Hyperplane_Arrangement}. That is, $P_W(\bm x) = \sum_{W\subseteq V_l} x_l,\ P_{\bm 0}(\bm x) = \sum_{i=1}^r x_i$.
	
	\vspace{0.5em}

	In \cite[Lemma 4.4-4.6]{MPVIHA}, a special tuple $\bm u \in \mathbb Q_{<0}^{r}$ was constructed. To adjust the notation there for this paper, we should substitute $d$, $n+1$, and $a_i-1$ in \cite{MPVIHA} with $r$, $n$, and $u_i$ here. $\bm u$ satisfies the following.
	
	\noindent(1) $\sum_{i=1}^r u_i = -n$, see the statement of \cite[Lemma 4.6]{MPVIHA}.
	
	\noindent(2) $b_W := \nu_j + P_W(\bm u) > 0$ for all $W \in \mathcal S\setminus \{\bm 0\}$, see the last line of \cite[Lemma 4.6]{MPVIHA}. 
	
	Then by (1), $\nu_{\bm 0} + P_{\bm 0}(\bm u) = n+\sum_{i=1}^r a_i = 0$. Take $m$ to be a positive integer such that $c_i := -mu_i \in \mathbb Z_{>0}$ for all $i$. Then
	\begin{displaymath}
		P_{W}(\bm c)\cdot (-\frac{\nu_{\bm 0}}{P_{\bm 0}(\bm c)}) + \nu_W = P_W(\bm u)+\nu_0 > 0,
	\end{displaymath}
	for all $W \in \mathcal S\setminus \{\bm 0\}$.

	\begin{proof}[Proof of Theorem \ref{Generically_n_over_d}]
		Since $\bm c$ is a good tuple at $\bm 0 \in \mathcal S$, Theorem \ref{Generically_n_over_d} follows immediately from Theorem \ref{generic existence}.
	\end{proof}

	Theorem \ref{Generically_n_over_d} leads to a generic result about the strong monodromy conjecture for hyperplane arrangements. That is, Corollary \ref{generic monodromy conjecture}.
	\begin{proof}[Proof of Corollary \ref{generic monodromy conjecture}]
		
		Using the same argument as Remark 2.5 of \cite{Monodromy_Conjecture_for_Hyperplane_Arrangement_Budur_Mustata_Teitler}, we reduce the problem to central hyperplane arrangement cases.
		
		Denote by $\mathcal S_{\mathrm{den}}$ the set of all dense edges of $A$. For $W\in \mathcal{S}_{\mathrm{den}}$, let $\mathcal P_W$ be the set of indices $l$ such that $W\subseteq \{f_l = 0\}$. 
		
		From the log resolution by blowing up all dense edges, we know that for $\bm a\in \mathbb N_{>0}^r$, all candidate poles of $Z_{\bm f^{\bm a}}^{\mathrm{mot}}(s)$ are of the form $\frac{\nu_W}{P_W(\bm a)}$, for some $W\in \mathcal S_{\mathrm{den}}$, see \cite[Proof of part (a) of Theorem 1.3]{Monodromy_Conjecture_for_Hyperplane_Arrangement_Budur_Mustata_Teitler}. Let $f_W := \prod_{l\in \mathcal P_W} f_l^{a_l}$, which defines a hyperplane arrangement in $\mathbb C^n/W$. Note that if we choose $x \in W \setminus \bigcup_{l\notin \mathcal P_W} V_l$, then we have an isomorphism between germs of log pairs $(\mathbb C^n,f,x) \simeq ((\mathbb C^{n}/W)\times \mathbb C^{n-\nu_W},f_W\otimes 1,(\bm 0',\bm 0''))$, where $\bm 0'$ and $\bm 0''$ are origins of $\mathbb C^n/W$ and $W$, respectively. Since $f_W$ is again homogeneous, by {Remark \ref{basic facts about B-fun}}, $b_{f,x}(s) = b_{f_W,\bm 0'}(s) = b_{f_W}(s)$. So $b_{f_W}(s)$ is a factor of $b_f(s)$.
		
		For $W\in \mathcal S_{\mathrm{den}}$, we have a biholomorphic map
		\begin{displaymath}
			\{\bm x \in \mathbb U^r \mid \sum_{i=1}^r x_i = 1\} \overset{\simeq}{\longrightarrow} \{\bm x \in \mathbb U^r \mid P_W(\bm x) = 1\}, \bm x \mapsto \frac{\bm x}{P_W(\bm x)}.
		\end{displaymath}
		Applying {Theorem \ref{Generically_n_over_d}} to $(\mathbb C^n/W, \{\prod_{l\in \mathcal P_W} f_l = 0\})$ and using the above isomorphism, there exists an analytically open subset $U_{W} \subseteq \{\bm x \in \mathbb U^r \mid \sum_{i=1}^r x_i = 1\}$, such that if $\bm a \in \mathbb N_{>0}^r$ and $(\frac{a_1}{\sum_{i=1}^r a_i},...,\frac{a_r}{\sum_{i=1}^r a_i}) \in U_W$, then $-\nu_W/P_W(\bm a)$ is a root of $b_{f_W}(s)$. Taking $U = \bigcap_{W\in\mathcal S_{\mathrm{den}}} U_W$, we conclude the proof.
	\end{proof}
	
	From our strategy for $n/d$-conjecture, it is natural to ask the following.
	\begin{question}
		Let $f_1,...,f_r\in \mathbb C[\bm z]\setminus \mathbb C$ be polynomials of degree one, and $f = f_1^{b_1}\cdots f_r^{b_r},b_i \in \mathbb N_{>0}$ for all $i$, be of degree $d$. Suppose $\{f = 0\}$ is a central essential indecomposable hyperplane arrangement. Does there always exist a $\varphi \in C_c^{\infty}(\mathbb C^n)$ such that $-n/d$ is a pole of $Z_{f,\varphi}(s)$?
	\end{question}
	
	In the $n=2$ case, we will give a positive answer to this problem in the next section through explicitly constructing one $\varphi$  and analyzing the residue at $-2/d$. However, the question in higher-dimensional cases turns out to be much more complicated. 
	


	\section{Dimension two case}\label{section4}
	
	Suppose $f_1,...,f_r \in \mathbb C[x,y]\setminus \mathbb C$ are linear functions and $\{f_1f_2\cdots f_r = 0\}$ defines a central essential indecomposable hyperplane arrangement. After a linear coordinate change and scaling, we may assume $f_i = x+a_iy$, where $a_1 = 0$ and $a_2,...,a_r$ are pairwise distinct nonzero complex numbers, $r \geq 2$. Let $\bm b\in \mathbb N_{>0}^r$, and assume that $b_1 = \max \{b_i\}_{i=1}^r$. 
	
	
	Let $\varphi(x,y) = \alpha(x)\beta(y)$, where $\alpha,\beta \in C_c^{\infty}(\mathbb C)$, $\alpha|_{\mathbb D_1(0)} = \beta|_{\mathbb D_1(0)} = 1 , \alpha,\beta \geq 0$. We point out that the radius $1$ is not essential, since $\bm f^{\bm b} = f_1^{b_1}\cdots f_r^{b_r}$ is homogeneous and $Z_{\bm f^{\bm b},\varphi}(s)$ only changes by the multiple $\lambda^{-\vert b\vert s-n}$ if one scales $\varphi$ by $\varphi_{\lambda}: (x,y) \mapsto \varphi(\lambda x,\lambda y),\lambda > 0$. It does not affect the survival of poles or the arguments of residues.
	
	\begin{proposition}\label{easier side}
		The order of $-2/d$ as a pole of $Z_{\bm f^{\bm b},\varphi}(s)$ is not greater than $2$ if $b_1 = d/2$ and is not greater than $1$ if $b_1\neq d/2$. Moreover, if $b_1 \leq d/2$ and $-2/d$ is not an order two pole of $Z_{\bm f^{\bm a}, \varphi}(s)$, then it is an order one pole with a positive residue.
	\end{proposition}
	\begin{proof}
		From the set of candidate poles, i.e. $\Lambda \cap (\{\bm b\}\times \mathbb C)$ in {Proposition \ref{meromorphic extensino of VAZF}}, we see that the order of $-2/d$ as a pole is $\leq 2$. Moreover, the order is $\leq 1$ if $b_1 < d/2$. The case for $b_1 > d/2$ can be covered by our computation of residues later. In the $b_1 > d/2$ case, we will show that the value of $(ds+2)Z_{\bm f^{\bm b},\varphi}(s)$ at $-2/d$ is finite and hence the order is $\leq 1$. 
		
		For the second assertion, the argument is based on {Remark \ref{LCT_LARGEST_POLE}}. Since $-2/d$ is the log canonical threshold if $b_1\leq d/2$, it is the largest pole of $Z_{\bm f^{\bm b},\varphi}(s)$.   If $-2/d$ is not an order two pole, then $Z_{\bm f^{\bm b},\varphi}(s)$ has a positive residue at $-2/d$ since $\displaystyle\lim_{\substack{s\in \mathbb R, s \to (-2/d)^+}} Z_{\bm f^{\bm b},\varphi}(s) \to +\infty$. 
	\end{proof}
	

	So, it suffices to only consider the case $b_1 > d/2$. From now on, we fix $b_1,d\in \mathbb N_{>0}$ with $b_1 < d < 2b_1$. For every $\bm b' := (b_2,...,b_r) \in \mathbb R_{>0}^{r-1}$, let $\zeta_{\bm b'}(-)$ denote the meromorphic extension of $Z_{\varphi}(b_1,\bm b',-) : \mathbb U \to \mathbb C$. We note that $Z_{\bm f^{\bm b},\varphi}(s) = \zeta_{\bm b'}(s)$ if $\bm b' \in \mathbb N_{>0}^{r-1}$.
	\begin{remark}
		The meromorphy of $\zeta_{\bm b'}(s)$ follows from the same proof as {Proposition \ref{meromorphic extensino of VAZF}}. 
	\end{remark}

	Our strategy is to consider the residue of $\zeta_{\bm b'}(s)$ at $-2/d$ when $\bm b' = (b_2,...,b_r)$ varies in $\overline{\mathbb G} \subseteq \mathbb R^{r-1}$, where $
	\mathbb G := \{\bm b' \in \mathbb R_{> 0}^{r-1} \mid \sum_{i=2}^r b_i = d-b_1\}$ is convex and $\overline{\mathbb G}$ is its closure.
	

	A log resolution of $f$ is obtained by blowing up the origin of $\mathbb C^2$. By the change of variables $x \mapsto uv, y \mapsto  v$, for $(\bm b',s) \in \mathbb G\times  \mathbb U$, we have
	\begin{align*}
		\zeta_{\bm b'}(s) & = (\frac{\sqrt{-1}}{2})^2\int_{\mathbb C^2} \vert v\vert^{2(ds+1)} \vert u\vert^{2b_1s} \prod_{l = 2}^r \vert u+a_l\vert^{2b_ls}\alpha(uv)\beta(v)\, \mathrm{d} u\wedge \mathrm{d} \bar u\wedge\mathrm{d} v\wedge \mathrm{d} \bar v.
	\end{align*}
	\begin{remark}
		Take $0 < \varepsilon \ll \big((d+1)(n+1)\big)^{-1}$ and set $\mathbb U_{\varepsilon} := \{ s \in \mathbb C \mid 0 < \mathrm{Re}\, s < \varepsilon\}$. In fact, the right-hand side is a holomorphic function of $(\bm b',s)$ on a neighborhood of $\overline{\mathbb G} \times \mathbb U$. It can be verified by the same dominated convergence argument as in {Proposition \ref{meromorphic extensino of VAZF}}.
	\end{remark}
	
	
	
	Let $v = \rho e^{\theta\sqrt{-1}}$, then
	\begin{displaymath}
		\zeta_{\bm b'}(s) = \int_{0}^{+\infty} \int_{0}^{2\pi}\rho^{2ds+3}\beta(v)\bigg( \frac{\sqrt{-1}}{2}\int_{\mathbb C} \vert u\vert^{2b_1s}\prod_{l = 2}^r \vert u+a_l\vert^{2b_ls}\alpha(uv)\,\mathrm{d} u\wedge \mathrm{d} \bar u \bigg) \, \mathrm{d} \rho\wedge \mathrm{d} \theta.
	\end{displaymath}
	For $s \in \mathbb U$, let
	\begin{displaymath}
		A(\bm b',s,v) := \frac{\sqrt{-1}}{2}\int_{\mathbb C}\vert u\vert^{2b_1s}\prod_{l = 2}^r \vert u+a_l\vert^{2b_ls}\alpha(uv)\mathrm{d} u\wedge \mathrm{d} \bar u.
	\end{displaymath}
	For fixed $v \neq 0$ and $\bm b' \in \mathbb G$, $A(\bm b',-,v)$ can also be meromorphically extended to the whole $\mathbb C$ again through the function equations obtained by integration by parts. Moreover, we will see in the proof of the upcoming {Lemma \ref{Properties_of_A1}} that $A(\bm b',-2/d,v)$ is finite for all $\bm b' \in \overline{\mathbb G}$.
	
	
	\vspace{0.5em}

	The ultimate goal of the rest of this section is to prove $\mathrm{Res}|_{s = -2/d}\, \zeta_{\bm b'}(s) \in \mathbb R_{<0}$. We sketch the proof before proving the technical lemmas. First, we characterize the behavior of $A(\bm b',-2/d,v)$ when $v \to 0$, see {Lemma \ref{Properties_of_A1}}. Second, we prove that the residue $\mathrm{Res}|_{s = -2/d}\, \zeta_{\bm b'}(s)$ is equal to $\displaystyle 2\pi \lim_{v\to 0} A(\bm b',-2/d,v)$, see {Lemma \ref{Properties_of_A2}} and {Lemma \ref{Properties_of_A3}}.  Third, we perform a hard analysis on $\displaystyle \lim_{v\to 0} A(\bm b',-2/d,v)$ and find it negative, see {Lemma \ref{Properties_of_A4}}. Finally, we conclude Lemma \ref{Properties_of_A1}-Lemma \ref{Properties_of_A4} with Proposition \ref{Summary of A1-A4}. Adding up Proposition \ref{Summary of A1-A4} and {Proposition \ref{easier side}}, we obtain {Theorem \ref{n_over_d_in_Dimension_two}}.

	\begin{lemma}\label{Properties_of_A1}
		For a fixed $\bm b'\in \overline{\mathbb G}$, $\displaystyle \lim_{v \to 0} A(\bm b',-2/d,v)$ exists. Moreover, it can be expressed as
		\begin{align*}
			& \lim_{v \to 0} A(\bm b',-2/d,v) \\
			& = \lim_{\delta \to 0} \bigg( \frac{\sqrt{-1}}{2}\int_{\vert u\vert > \delta}\vert u\vert^{-4b_1/d}\prod_{l = 2}^r \vert u+a_l\vert^{-4b_l/d}\mathrm{d} u\wedge \mathrm{d} \bar u + \frac{2\pi\delta^{-4b_1/d+2}}{-4b_1/d+2} \prod_{l = 2}^{r} \vert a_l\vert^{-4b_l/d} \bigg).
		\end{align*}
	\end{lemma}

	Let $\mathbb W := \{s\in \mathbb C\mid \mathrm{Re}\, s > -2/d-\varepsilon,s\neq -1/b_1 \}, 0 < \varepsilon \ll \big((d+1)(n+1)\big)^{-1}$.
	\begin{lemma}\label{Properties_of_A2}
		Fix $\bm b' \in \overline{\mathbb G}$ and we write $v = \rho e^{\theta\sqrt{-1}}$ in the form of polar coordinates.
		
		\noindent(1) For fixed $s\in \mathbb W$, $A(\bm b',s,-)$ is differentiable on $\mathbb C^*$. 
		
		\noindent(2) For fixed $v \neq 0$, $\frac{\partial}{\partial \rho} A(\bm b',-,v)$ is holomorphic on $\mathbb W$. 
		
		\noindent(3) $\frac{\partial}{\partial \rho}A(\bm b',-,-),\frac{\partial}{\partial s}\frac{\partial}{\partial \rho} A(\bm b',-,-)$ are continuous on $\mathbb W\times \mathbb C^*$. 
		
		\noindent(4) The following five functions are integrable in any bounded neighborhood of $0 \in \mathbb C$. 
		\begin{align*}
			& v\mapsto \vert v\vert^{2ds+4}\frac{\partial}{\partial \rho} A(\bm b',s,v),\ v\mapsto \vert v\vert^{2ds+4}\frac{\partial}{\partial s}\frac{\partial}{\partial \rho} A(\bm b',s,v),\\
			& v\mapsto \vert v\vert^{2ds+4} A(\bm b',s,v),\ v\mapsto \vert v\vert^{2ds+4}(\log \vert v\vert) \frac{\partial}{\partial s}\frac{\partial}{\partial \rho} A(\bm b',s,v),\ v\mapsto \vert v\vert^{2ds+4}\frac{\partial}{\partial s} A(\bm b',s,v).
		\end{align*}
		Moreover, for fixed $s_0\in \mathbb W$, there are a neighborhood $U_{s_0}$ and a function $g$ integrable in any bounded neighborhood of $0$ such that the above five functions are all dominated by $\vert g\vert$ when $s \in U_{s_0}$.
		
		%
	\end{lemma}
	\begin{lemma}\label{Properties_of_A3}
		Fix $\bm b' \in \overline{\mathbb G}$. For $s\in \mathbb U$, we can do integration by parts to $\rho$ in the integral expression of $\zeta_{\bm b'}(s)$. Thus, by {Lemma \ref{Properties_of_A2}}, for $s \in \mathbb W$, we have
		\begin{align*}
			\zeta_{\bm b'}(s) = -\frac{1}{2ds+4}\int_{0}^{+\infty} \int_{0}^{2\pi}\rho^{2ds+4} \frac{\partial}{\partial \rho}(\beta(v)A(\bm b',s,v))\, \mathrm{d}\rho \mathrm{d} \theta,
		\end{align*} 
		where $v = \rho e^{\theta\sqrt{-1}}$, where $\rho \geq 0,\theta \in [0,2\pi)$.  
		
		Consequently, the residue of $\zeta_{\bm b'}(s)$ at $-2/d$ is given by
		\begin{displaymath}
			\mathrm{Res}|_{s = -2/d}\, \zeta_{\bm b'}(s) = -\int_{0}^{+\infty} \int_{0}^{2\pi} \frac{\partial}{\partial \rho}(\beta(v)A(\bm b',-2/d,v))\, \mathrm{d}\rho \mathrm{d} \theta = 2\pi \lim_{v \to 0} A(\bm b',-2/d,v).
		\end{displaymath}
	\end{lemma}
	\begin{remark}
		$2\pi \lim_{v \to 0} A(\bm b',-2/d,v)$ is a real-valued function of $\bm b'$.
	\end{remark}
	\begin{remark}
		{Lemma \ref{Properties_of_A2}} is a preparatory lemma for {Lemma \ref{Properties_of_A3}}.
	\end{remark}

	\begin{lemma}\label{Properties_of_A4}
		For all $\bm b' \in \overline{\mathbb G}$, $\lim_{v \to 0} A(\bm b',-2/d,v) < 0$.
	\end{lemma}

	\begin{proposition}\label{Summary of A1-A4}
		The residue of $\zeta_{\bm b'}(s)$ at $-2/d$ is a negative real number.
	\end{proposition}
	So far, all the cases in dimension two have been settled and thus Theorem \ref{n_over_d_in_Dimension_two} is proved.

	Now we start to prove the lemmas.

	\begin{proof}[Proof of {Lemma \ref{Properties_of_A1}}]
		For $\mathrm{Re}\, s > -d/2-\varepsilon\ ,0<\varepsilon \ll 1$, we can compute $A(\bm b',s,v)$ as follows. Taking $0 < \delta \ll \min \{\vert b_i\vert/2\}_{i=2}^r$, we see that
		
		\begin{displaymath}
			A_1(\bm b',s,v,\delta) := \frac{\sqrt{-1}}{2}\int_{\vert u\vert > \delta}\vert u\vert^{2b_1s}\prod_{l = 2}^r \vert u+a_l\vert^{2b_ls}\alpha(uv)\mathrm{d} u\wedge \mathrm{d} \bar u.
		\end{displaymath}
		is holomorphic for $\mathrm{Re}\, s > -\min\{{1}/{b_l}\}_{l\geq 2}$ by the same dominated convergence as in {Proposition \ref{meromorphic extensino of VAZF}}. Since $-\min\{{1}/{b_l}\}_{l\geq 2} < -2/d-\varepsilon$, $A_1(\bm b',-,v,\delta)$ is holomorphic on $\mathbb W$. Moreover, we have $\vert u\vert^{2b_1(-2/d)}\prod_{l = 2}^r \vert u+a_l\vert^{2b_l(-2/d)} \in L^1(\mathbb C\setminus \mathbb D_{\delta}(0))$ since $b_l(-2/d) > -2$ for all $l=2,...,r$ and $-4\sum_{i=1}^r b_i/d = -4$. Then, taking $s = -2/d$, since $\alpha$ is bounded, by the dominated convergence, we have
		\begin{displaymath}
			\lim_{v \to 0} A_1(\bm b',-2/d,v) = \frac{\sqrt{-1}}{2}\int_{\vert u\vert > \delta}\vert u\vert^{-4b_1/d}\prod_{l = 2}^r \vert u+a_l\vert^{-4b_l/d}\mathrm{d} u\wedge \mathrm{d} \bar u.
		\end{displaymath}
		
		So we only have to consider the meromorphic extension of
		\begin{displaymath}
			A_2(\bm b',s,v,\delta) := \frac{\sqrt{-1}}{2}\int_{\vert u\vert < \delta}\vert u\vert^{2b_1s}\prod_{l = 2}^r \vert u+a_l\vert^{2b_ls}\alpha(uv)\mathrm{d} u\wedge \mathrm{d} \bar u.
		\end{displaymath}
		Let $u = \lambda e^{\gamma\sqrt{-1}}$, where $\lambda \geq 0, \gamma \in [0,2\pi)$. Integrating by parts when $\mathrm{Re}\, s > 0$, we have
		\begin{align*}
			& A_2(\bm b',s,v,\delta)  = \int_0^{2\pi} \int_{0}^\delta \lambda^{2b_1s+1} \prod_{l = 2}^r \vert u+a_l\vert^{2b_ls} \alpha(uv)\mathrm{d} \lambda\mathrm{d}  \gamma \\
			& = \frac{\delta^{2b_1s+2}}{2b_1s+2}\int_{0}^{2\pi} \prod_{l = 2}^r \vert \delta e^{\gamma\sqrt{-1}}+a_l\vert^{2b_ls} \alpha(\delta e^{\gamma\sqrt{-1}}v) \mathrm{d}  \gamma \\
			& \quad -\frac{1}{2b_1s+2} \int_0^{2\pi} \int_{0}^\delta\lambda^{2b_1s+2} \frac{\partial}{\partial \lambda} \bigg(\prod_{l = 2}^r \vert \lambda e^{\gamma\sqrt{-1}}+a_l\vert^{2b_ls} \alpha(\lambda e^{\gamma\sqrt{-1}}v) \bigg)\, \mathrm{d} \lambda\mathrm{d}  \gamma. \tag{\ref{Properties_of_A1}.1}
		\end{align*}
		Let $A_{21}(\bm b',s,v,\delta)$ and $A_{22}(\bm b',s,v,\delta)$ be the two terms on the right-hand side of (\ref{Properties_of_A1}.1), respectively. We show that $\displaystyle \lim_{v \to 0} A_2(\bm b',-2/d,v)$ is the evaluation of the right-hand side at $s = -2/d$ and $v = 0$. It suffices to prove this claim for $A_{21}$ and $A_{22}$.
		
		The assertion for $A_{21}$ is easy since everything is bounded. For $A_{22}$, we first notice that in the Taylor expansion of the following function at $0$, there is no degree-one term.
		$$\lambda \mapsto \int_0^{2\pi} \prod_{l = 2}^r \vert \lambda e^{\gamma\sqrt{-1}}+a_l\vert^{2b_l(-2/d)} \alpha(\lambda e^{\gamma\sqrt{-1}}v) \, \mathrm{d}  \gamma.$$
		Indeed, since $\vert 1+a \lambda\vert^b = 1 + b(a+\bar a)\lambda + O(\lambda^2)$ and $ 2b_l(-2/d)a_l^{-1}\int_{0}^{2\pi} \cos \gamma\, \mathrm{d} \gamma = 0$, we see that $\vert \lambda e^{\gamma\sqrt{-1}}+a_l\vert$ contributes no degree one term. In addition, $\alpha(\lambda e^{\gamma\sqrt{-1}})$ is constant when $\lambda \in [0,1]$. Hence, 
		\begin{displaymath}
			\frac{\partial}{\partial \lambda}\int_0^{2\pi} \prod_{l = 2}^r \vert \lambda e^{\gamma\sqrt{-1}}+a_l\vert^{2b_l(-2/d)} \alpha(\lambda e^{\gamma\sqrt{-1}}v) \, \mathrm{d}  \gamma = \lambda \cdot \mathcal E(\lambda,v),
		\end{displaymath}
		for some smooth function $\mathcal E$ near $(0,0)$. Since $2b_1(-2/d)+3 > -1$, $\lambda^{2b_1(-2/d)+3}$ is integrable on the interval $(0,\delta)$. So, we have
		\begin{displaymath}
			A_{22}(\bm b',-2/d,v,\delta) = -\frac{1}{2b_1(-2/d)+2}\int_{0}^\delta \lambda^{2b_1(-2/d)+3} \mathcal E(\lambda,v) \, \mathrm{d} \lambda.
		\end{displaymath}
		Since $\mathcal E$ is continuous near $(0,0)$, there exists an upper bound $M$ of $\mathcal E$ on $[0,\delta] \times \mathbb D_{\varepsilon '}(0)$ for small $\varepsilon' > 0$. Hence,
		\begin{displaymath}
			\limsup_{v \to 0}\vert \int_{0}^\delta \lambda^{2b_1(-2/d)+3} \mathcal E(\lambda,v) \, \mathrm{d} \lambda\vert \leq \frac{M}{2b_1(-2/d)+3} \delta^{2b_1(-2/d)+4} \to 0, \tag{\ref{Properties_of_A2}.2}
		\end{displaymath}
		as $\delta \to 0$. In conclusion, only $A_1$ and $A_{21}$ survive when $\delta \to 0$, i.e.	
		
		%
		\begin{align*}
			& \lim_{v \to 0} A(\bm b',-2/d,v) \\
			& = \lim_{\delta \to 0} \bigg( \frac{\sqrt{-1}}{2}\int_{\vert u\vert > \delta}\vert 	u\vert^{-4b_1/d}\prod_{l = 2}^r \vert u+a_l\vert^{-4b_l/d}\mathrm{d} u\wedge \mathrm{d} \bar u + \frac{2\pi\delta^{-4b_1/d+2}}{-4b_1/d+2} \prod_{l = 2}^{r} \vert a_l\vert^{-4b_l/d} \bigg).
		\end{align*}
	\end{proof}
	\begin{proof}[Proof of {Lemma \ref{Properties_of_A2}}]
		The proof is a refinement of the procedure of {Lemma \ref{Properties_of_A1}}. The notation is as in Lemma \ref{Properties_of_A1} and we have
		\begin{align*}
			A_1(\bm b',s,v,\delta) & = \frac{\sqrt{-1}}{2}\int_{\vert u\vert > \delta}\vert u\vert^{2b_1s}\prod_{l = 2}^r \vert u+a_l\vert^{2b_ls}\alpha(uv)\mathrm{d} u\wedge \mathrm{d} \bar u,\\
			A_2(\bm b',s,v,\delta) & = \frac{\delta^{2b_1s+2}}{2b_1s+2}\int_{0}^{2\pi} \prod_{l = 2}^r \vert \delta e^{\gamma\sqrt{-1}}+a_l\vert^{2b_ls} \alpha(\delta e^{\gamma\sqrt{-1}}v) \mathrm{d}  \gamma \\
			& \quad -\frac{1}{2b_1s+2} \int_0^{2\pi} \int_{0}^\delta\lambda^{2b_1s+2} \frac{\partial}{\partial \lambda} \bigg(\prod_{l = 2}^r \vert \lambda e^{\gamma\sqrt{-1}}+a_l\vert^{2b_ls} \alpha(\lambda e^{\gamma\sqrt{-1}}v) \bigg)\mathrm{d} \lambda\mathrm{d}  \gamma \\
			& = A_{21}+A_{22},
		\end{align*}
		where $u = \lambda e^{\gamma\sqrt{-1}}$. Then (1)(2)(3) follows from $\alpha'\in C_c^{\infty}(\mathbb C)$ and the dominated convergence. $\frac{\partial}{\partial \rho},\frac{\partial}{\partial s}$ hence commutes with the integral symbol.

		For (4), we fix $0 < \delta \ll 1$ and $s_0 \in \mathbb W$. Using an estimate similar to (\ref{Properties_of_A2}.2), all integrability and dominating with respect to $A_{22}$ can be verified. The case is much easier for the $A_{21}$ because all $\vert \lambda e^{\gamma\sqrt{-1}}+a_l\vert$ with $l=2,...,r$ are positively bounded from below. So, it suffices to prove the domination for $A_1$.

		We only prove the domination for $ v\mapsto \vert v\vert^{2ds+4}\frac{\partial}{\partial \rho} A(\bm b',s,v), v = \rho e^{\theta \sqrt{-1}}$. The other four functions can be dominated by using the same method. Choose a neighborhoods $U_{s_0} \subseteq \mathbb W$ of $s_0$ such that $\overline {U}_{s_0} \subseteq \mathbb W$, then there exists $M,\delta' > 1$ such that $\vert \frac{\partial A_1}{\partial \rho}(\bm b',s',v,\delta)\vert \leq M\big(1+\vert \int_{\delta ' < \vert u\vert < w\vert v\vert ^{-1}} \vert u\vert^{2\mathrm{Re}(ds')+1}\,\mathrm{d} u\wedge \mathrm{d} \bar u\vert \big)$ for all $s'\in U_{s_0}$. Here $w > 1$ is a constant such that $\mathrm{Supp}\, \beta \subset \mathbb D_{w}(0)$ . The ``$+1$'' on the exponent comes from the estimation $\vert\frac{\partial}{\partial  \rho} \big(\alpha(u\rho e^{\theta\sqrt{-1}})\big) \vert = \vert ue^{\theta\sqrt{-1}}\frac{\partial \alpha }{\partial z} + \bar u e^{-\theta\sqrt{-1}} \frac{\partial \alpha}{\partial \bar z}\vert \leq \big( \Vert \frac{\partial\alpha}{\partial z}\Vert_{L^{\infty}}+ \Vert \frac{\partial\alpha}{\partial \bar z}\Vert_{L^{\infty}}\big) \vert u\vert $. Consequently, we have $\vert \frac{\partial A_1}{\partial \rho}(\bm b',s',v,\delta)\vert \leq M\big(1+ w\vert v\vert^{-1}(1+(w\vert v\vert^{-1})^{2\mathrm{Re}(ds')+1})\big)$. So $\vert v\vert^{2ds'+4}\vert \frac{\partial A_1}{\partial \rho}(\bm b',s',v,\delta)\vert$ is dominated by $g(s',v) := M \vert v\vert^{2\mathrm{Re}(ds')+4} + Mw \vert v \vert^{2\mathrm{Re}\,(ds')+3} + Mw^{2\mathrm{Re}(ds')+2}\vert v \vert^{2}$. Let $\epsilon_0 = \inf_{s' \in U_{s_0}} 2\mathrm{Re}(ds')+4$ and $\epsilon_1 = \sup_{s'\in U_{s_0}} 2\mathrm{Re}(ds')+4$, then $g(s',v) \leq M(\vert v\vert^{\epsilon_0}+1)+Mw(\vert v\vert^{\epsilon_0-1}+1) + Mw^{\epsilon_1-2}\vert v\vert^2 =: g(v)$. Since $s'\in \mathbb W$, we have $\epsilon_0 > 2d\varepsilon > -1$. Then $g$ in integrable and hence is a dominating function. 
		%
	\end{proof}
	
	\begin{proof}[Proof of {Lemma \ref{Properties_of_A3}}]
		For $\mathrm{Re}\, s > 0$, by (4) of {Lemma \ref{Properties_of_A2}}, we can do integration by parts to $\rho$. Again by {Lemma \ref{Properties_of_A2}} (4) and since $\beta \in C_c^{\infty}(\mathbb C)$, we see that the map
		\begin{displaymath}
			s\mapsto \int_{0}^{+\infty} \int_{0}^{2\pi}\rho^{2ds+4} \frac{\partial}{\partial 	\rho}(\beta(v)A(\bm b',s,v))\, \mathrm{d}\rho  \mathrm{d} \theta
		\end{displaymath}
		is holomorphic on $\mathbb W$ for fixed $\bm b'$. We are done.
	\end{proof}
	
	Now we prove {Lemma \ref{Properties_of_A4}}. Basically, we are going to prove that $\lim_{v \to 0} A(\bm b',-2/d,v)$ is a strict convex function in a neighborhood of $\overline{\mathbb G}$ and hence reaches its maximum at the vertices of $\overline{\mathbb G}$. Then we prove at each vertex $\bm b'$ of $\overline{\mathbb G}$, $\lim_{v \to 0} A(\bm b',-2/d,v) = 0$.
	\begin{proof}[Proof of {Lemma \ref{Properties_of_A4}}]
		Let 
		\begin{align*}
			\tilde C(\bm b') := \lim_{\delta \to 0} \bigg( \frac{\sqrt{-1}}{2}\int_{\vert u\vert > \delta}\vert 	u\vert^{-4b_1/d}\prod_{l = 2}^r \vert u+a_l\vert^{-4b_l/d}\mathrm{d} u\wedge \mathrm{d} \bar u + \frac{2\pi\delta^{-4b_1/d+2}}{-4b_1/d+2} \prod_{l = 2}^{r} \vert a_l\vert^{-4b_l/d} \bigg).
		\end{align*}
		We normalize $\tilde C(\bm b')$ as $C(\bm b') := \prod_{l = 2}^r \vert a_l\vert^{4b_l/d} \tilde C(\bm b')$. It suffices to show $C(\bm b') < 0$ for all $\bm b'\in \mathbb G$. 
		\begin{displaymath}
			C(\bm b') = \lim_{\delta \to 0} \bigg( \frac{\sqrt{-1}}{2}\int_{\vert u\vert > \delta}\vert 	u\vert^{-4b_1/d}\prod_{l = 2}^r \vert a_l^{-1}u+1\vert^{-4b_l/d}\mathrm{d} u\wedge \mathrm{d} \bar u + \frac{2\pi\delta^{-4b_1/d+2}}{-4b_1/d+2} \bigg).
		\end{displaymath}
		Note that the right-hand side is a function on $\mathbb G_{\epsilon} = \{\bm b'\in \mathbb R_{>-\epsilon}^{r-1} \mid \vert \bm b'\vert = d-b_1 < \frac{d}{2}\}$ for $\epsilon$ small. This can be verified by noticing: (1) $A_1(\bm b',-2/d,v,\delta),A_2(\bm b',-2/d,\delta)$ in the proof of {Proposition \ref{Properties_of_A1}} are well defined on $\mathbb G_{\epsilon}$, (2) their sum is independent of $\delta$ (for $\delta_1<\delta_2$, use integration by parts on the annulus region), and (3) $A_{22}$ again tends to zero as $\delta \to 0$. 
		
		\noindent\textit{Claim 1: $\frac{\partial C}{\partial b_j}(\bm b')$ and $\frac{\partial^2 C}{\partial b_j\partial b_k}(\bm b')$ exist on $\mathbb G_\epsilon$.}
		
		For $\bm b'$ and $\tilde{\bm b}' = \bm b'+(\tilde b_2,0,...,0),0 < \vert b_l\vert < \epsilon \ll 1$ in $\mathbb G_{\epsilon}$, we have
		\begin{displaymath}
			C(\tilde{\bm b}') - C(\bm b') = \lim_{\delta \to 0}\frac{\sqrt{-1}}{2}\int_{\vert u\vert > \delta}\vert 	u\vert^{-4b_1/d}(\vert a_2^{-1}u+1\vert^{-4\tilde b_2/d}-1)\prod_{l = 2}^r \vert a_l^{-1}u+1\vert^{-4b_l/d}\mathrm{d} u\wedge \mathrm{d} \bar u.
		\end{displaymath} 
		Let $u = \lambda e^{\gamma \sqrt{-1}}$, then
		\begin{align*}
			RHS & = \lim_{\delta \to 0} \int_{\delta}^{\infty} \lambda^{-4b_1/d+1} \bigg(\int_{0}^{2\pi} (\vert a_2^{-1}\lambda e^{\gamma\sqrt{-1}}+1\vert^{-4\tilde b_2/d}-1)\prod_{l = 2}^r \vert a_l^{-1}\lambda e^{\gamma\sqrt{-1}}+1\vert^{-4b_l/d}\, \mathrm{d} \gamma\bigg)\,\mathrm{d} \lambda.
		\end{align*}
		Look at the Taylor expansion of $\displaystyle{\int_{0}^{2\pi} (\vert a_2^{-1}\lambda e^{\gamma\sqrt{-1}}+1\vert^{-4\tilde b_2/d}-1)\prod_{l = 2}^r \vert a_l^{-1}\lambda e^{\gamma\sqrt{-1}}+1\vert^{-4b_l/d}\, \mathrm{d} \gamma}$ at $\lambda = 0$. The constant term is zero since $\big(\vert a_2^{-1}\lambda e^{\gamma\sqrt{-1}}+1\vert^{-4\tilde b_2/d}-1\big)|_{\lambda = 0} = 0$. The order one term is also zero since $\int_{0}^{2\pi} e^{\gamma\sqrt{-1}}\, \mathrm{d} \gamma = 0$. Therefore, around $0$, the function inside $\int_{\delta}^{\infty}$ is approximately some constant multiple of $\lambda^{-b_1/d+3}$, which is integrable near $0$. Moreover, since $\sum_{l=1}^r -4b_l/d -4\tilde b_2/d = -4-\tilde b_2/d < -3$ and $-4b_l/d -4\tilde b_2/d > -2$, this function is also integrable near each $-a_l$ and $\infty$. Hence, 
		\begin{displaymath}
			C(\tilde{\bm b}') - C(\bm b') = \int_{0}^{\infty} \lambda^{-4b_1/d+1} \bigg(\int_{0}^{2\pi} (\vert a_2^{-1}\lambda e^{\gamma\sqrt{-1}}+1\vert^{-4\tilde b_2/d}-1)\prod_{l = 2}^r \vert a_l^{-1}\lambda e^{\gamma\sqrt{-1}}+1\vert^{-4b_l/d}\, \mathrm{d} \gamma\bigg)\,\mathrm{d} \lambda.
		\end{displaymath}
		Then we consider ${(C(\tilde{\bm b}') - C(\bm b'))}/{\tilde b_2}$. Let $F(\lambda,\tilde b_2)$ denote the function inside $\int_0^{\infty}$ in the expression of ${(C(\tilde{\bm b}') - C(\bm b'))}/{\tilde b_2}$, i.e. ${(C(\tilde{\bm b}') - C(\bm b'))}/{\tilde b_2} = \int_0^\infty F(\lambda,\tilde b_2)\,\mathrm{d}\lambda$. We first show 
		\begin{displaymath}
			\int_{0}^{2\pi} \big({\vert a_2^{-1}\lambda e^{\gamma\sqrt{-1}}+1\vert^{-4\tilde b_2/d}-1}\big)\, \mathrm{d} \gamma = O(\lambda^2\tilde b_2).\tag{\ref{Properties_of_A4}.1}
		\end{displaymath}
		Taking $p = -4\tilde b_2/d$, we have
		\begin{displaymath}
			LHS = -p \int_0^{1}  \big(\int_0^{2\pi} ({\vert a_2^{-1}\lambda e^{\gamma\sqrt{-1}}+1\vert^{pt}}) \log \vert a_2^{-1}\lambda e^{\gamma\sqrt{-1}}+1\vert \, \mathrm{d} \gamma \big)\mathrm{d} t.
		\end{displaymath}
		Again, by the Taylor expansion, together with $\log \vert 1+u\vert = \frac{u+\bar u}{2} + O(\vert u\vert^2)$ and $\int_{0}^{2\pi} e^{\gamma \sqrt{-1}}\,\mathrm{d} \gamma = 0$, we find $\int_0^{2\pi} ({\vert a_2^{-1}\lambda e^{\gamma\sqrt{-1}}+1\vert^{pt}}) \log \vert a_2^{-1}\lambda e^{\gamma\sqrt{-1}}+1\vert \, \mathrm{d} \gamma = O(\lambda^2)$. (\ref{Properties_of_A4}.1) then follows from the continuity of $(\lambda,t) \mapsto \frac{1}{\lambda^2} \int_0^{2\pi} ({\vert a_2^{-1}\lambda e^{\gamma\sqrt{-1}}+1\vert^{pt}}) \log \vert a_2^{-1}\lambda e^{\gamma\sqrt{-1}}+1\vert \, \mathrm{d} \gamma$. Therefore, 
		
		\begin{displaymath}
			\frac{1}{\tilde b_2} \int_{0}^{2\pi} (\vert a_2^{-1}\lambda e^{\gamma\sqrt{-1}}+1\vert^{-4\tilde b_2/d}-1)\prod_{l = 2}^r \vert a_l^{-1}\lambda e^{\gamma\sqrt{-1}}+1\vert^{-4b_l/d}\, \mathrm{d} \gamma
		\end{displaymath}
		can be written in the form of $\lambda^2 \mathcal E(\lambda,\tilde b_2)$, with $\mathcal E$ continuous in a neighborhood $[0,2\delta_0) \times (-2\delta_0,2\delta_0),0<\delta_0\ll 1$ of $(0,\tilde b_2)$. Let $\mathcal E_0(\lambda) := \sup_{x \in [-\delta_0,\delta_0]}\vert \mathcal E(\lambda,x)\vert$, which is bounded on $[0,\delta_0]$. Moreover, $\lambda^{-4b_1/d+3}$ locally integrable near $0$, so $\lambda^{-4b_1/d+3}\mathcal E_0(\lambda)$ is a dominating function for $F\cdot {\bm 1}_{\{\lambda \leq \delta_0\}}$.
		
		For $\lambda > \delta_0$, we dominate $F(-,\tilde b_2)$ by the Lagrange inequality. $\vert \xi \vert < \epsilon$, we have
		\begin{displaymath}
			\vert a_2^{-1}\lambda e^{\gamma\sqrt{-1}}+1\vert^{-4\xi/d} < \vert a_2^{-1}\lambda e^{\gamma\sqrt{-1}}+1\vert^{-4\epsilon/d} + \vert a_2^{-1}\lambda e^{\gamma\sqrt{-1}}+1\vert^{4\epsilon/d}.
		\end{displaymath}
		By the Lagrange inequality, we see that $F\cdot {\bm 1}_{\{\lambda > \delta_0\}}$ is dominated by
		
		\begin{align*}
			& \lambda^{-4b_1/d+1} \bigg(\int_{0}^{2\pi} (\vert a_2^{-1}\lambda e^{\gamma\sqrt{-1}}+1\vert^{-4\epsilon/d} + \vert a_2^{-1}\lambda e^{\gamma\sqrt{-1}}+1\vert^{4\epsilon/d})\cdot \bigg\vert \log \vert a_2^{-1}\lambda e^{\gamma\sqrt{-1}}+1\vert\bigg\vert \\
			& \cdot \prod_{l = 2}^r \vert a_l^{-1}\lambda e^{\gamma\sqrt{-1}}+1\vert^{-4b_l/d}\, \mathrm{d} \gamma \bigg),
		\end{align*}
		which is integrable on $[\delta_0,+\infty)$.
		
		Hence, again by the dominated convergence, we have
		\begin{displaymath}
			\frac{\partial C}{\partial b_j} (\bm b') = \frac{\sqrt{-1}}{2}\int_{\mathbb C}\vert 	u\vert^{-4b_1/d}\cdot \log \vert a_j^{-1} u+1\vert\cdot  \prod_{l = 2}^r \vert a_l^{-1}u+1\vert^{-4b_l/d}\mathrm{d} u\wedge \mathrm{d} \bar u.
		\end{displaymath}
		By the same argument, we have the existence and formula for the second derivation.
		\begin{displaymath}
			\frac{\partial^2 C}{\partial b_j\partial b_k}(\bm b') = \frac{\sqrt{-1}}{2}\int_{\mathbb C}\vert 	u\vert^{-4b_1/d}\cdot \log \vert a_j^{-1} u+1\vert\cdot  \log \vert a_k^{-1} u+1\vert\cdot\prod_{l = 2}^r \vert a_l^{-1}u+1\vert^{-4b_l/d}\mathrm{d} u\wedge \mathrm{d} \bar u.
		\end{displaymath}
		
		\noindent\textit{Claim 2: $C$ is a convex function on $\mathbb G_{\epsilon}$.}
		
		Let $\bm X_u := (\log \vert a_2^{-1} u+1\vert,...,\log \vert a_r^{-1} u+1\vert)^t$, then the Hessian of $C$ is 
		$$\int_{\mathbb C} \vert u\vert^{-4b_1/d}\cdot ({\bm X}_u \cdot {\bm X}_u^t)\prod_{l = 2}^r \vert a_l^{-1}u+1\vert^{-4b_l/d}\,\frac{\sqrt{-1}}{2}\mathrm{d} u\wedge \mathrm{d} \bar u$$
		is an integral of positive semi-definite matrix. Hence, $\mathrm{Hess}(C)(\bm b')$ is a positive semi-definite matrix. Furthermore, suppose that there exist $\bm b'\in \mathbb G_{\epsilon}$ and $\bm w\in \mathbb R^{r-1}$ such that $\mathrm{Hess}(C)(\bm b')(\bm w,\bm w) = 0$, then $\bm X_u^t \cdot \bm w = 0$ for all $u \in \mathbb C\setminus \{0,-a_1,...,-a_r\}$. However, $\log \vert a_2^{-1} u+1\vert,...,\log \vert a_r^{-1} u+1\vert$ are $\mathbb C$-linearly independent,  and hence $\bm w = \bm 0$. Thus, we conclude that $\mathrm{Hess}(C)(\bm b')$ is positive definite for all $\bm b' \in \mathbb G_{\epsilon}$.
		
		Consequently, restrict $C$ to $\overline{\mathbb G}$, then the maximum of $C$ is reached at the vertices of $\overline{\mathbb G}$. Since $\overline{\mathbb G}$ is a simplex, its vertices are exact those $(b_2,...,b_r) \in \mathbb N^{r-1}$ such that $b_l = \begin{cases}d-b_1, & l =m,\\ 0, & l \neq m,\end{cases}$ for some $m=2,...,r$.
		
		\noindent\textit{Claim 3: For each vertex $\bm b_0'$ of $\overline{\mathbb G}$, $C(\bm b_0') = 0$ and thus the proof is finished.}
		
		Without loss of generality, we consider the vertex $\bm b_0' = (d-b_1,0,...,0) \in \mathbb N^{r-1}$. Then by {Lemma \ref{Properties_of_A3}} and $\zeta_{\bm b'_0}(s) = Z_{f_1^{b_1}f_2^{d-b_1},\varphi}(s)$,  we have
		\begin{displaymath}
			C(\bm b_0') = -\mathrm{Res}|_{s=-2/d}\, \zeta_{\bm b_0'}(s) = -\mathrm{Res}|_{s=-2/d}\, Z_{f_1^{b_1}f_2^{d-b_1},\varphi}(s).
		\end{displaymath}
		However, $-2/d$ is not a root of the $b$-function of $f_1^{b_1}f_2^{d-b_1}$, which is equal to $\prod_{i=1}^{b_1}(s+\frac{i}{b_1}) \cdot \prod_{j=1}^{d-b_1} (s+\frac{j}{d-b_1})$. Then it is not a pole by {Lemma \ref{Fundamental_Lemma_of_Archi_B_Function}} and the residue is thus zero.
	\end{proof}


\begin{thebibliography}{AAAAA00}
		
		\bibitem[ABCN$^+$17]{ABCN17}
		R. Artal Bartolo, P. Cassou-Nogu\`es, I. Luengo, and A. Melle-Hern\'andez.
		\newblock {Yano's conjecture for two-{P}uiseux-pair irreducible plane curve singularities}
		\newblock {\em Publ. Res. Inst. Math. Sci.}, 52(1): 211--239, 2017.
		
		\bibitem[AKMW02]{Torification_and_Factorization_of_Birational_Maps_Abramovich_Karu_Matsuki}
		D. Abramovich, K. Karu, K. Matsuki, and J. ~W{\l}odarczyk.
		\newblock Torification and factorization of birational maps.
		\newblock {\em J. Amer. Math. Soc.}, 15(3): 531--572, 2002.
		
		\bibitem[Ati70]{Continuation_Problem_of_Zeta_Function_Atiyah}
		M.~F. Atiyah.
		\newblock Resolution of singularities and division of distributions.
		\newblock {\em Comm. Pure Appl. Math.}, 23: 145--150, 1970.
		
		\bibitem[Bar84]{Bar84}
		D. Barlet.
		\newblock Contribution effective de la monodromie aux d\'eveloppements
		asymptotiques.
		\newblock {\em Ann. Sci. \'Ecole Norm. Sup. (4)}, 17(1): 293--315, 1984.
		
		\bibitem[Bar86]{Bar86}
		D. Barlet.
		\newblock Monodromie et p\^oles du prolongement m\'eromorphe de
		{$\int_X|f|^{2\lambda}\square$}.
		\newblock {\em Bull. Soc. Math. France}, 114(3): 247--269, 1986.
		
		
		\bibitem[Bat20]{Tame_Arrangement}
		D. Bath.
		\newblock Combinatorially determined zeroes of {B}ernstein-{S}ato ideals for
		tame and free arrangements.
		\newblock {\em J. Singul.}, 20: 165--204, 2020.
		
		
		\bibitem[Ber72]{Bernstein-polynomial_Continuations}
		I.~N. Bern\v{s}te\u{\i}n.
		\newblock Analytic continuation of generalized functions with respect to a
		parameter.
		\newblock {\em Funkcional. Anal. i Prilo\v zen.}, 6(4): 26--40, 1972.
		
		\bibitem[BG69]{Continuation_Problem_of_Zeta_Function_Bernstein_and_Gelfand}
		J.~N. Bern\v{s}te\v{\i}n and S.~I. Gel'fand.
		\newblock Meromorphy of the function {$P\sp{\lambda }$}.
		\newblock {\em Funkcional. Anal. i Prilo\v{z}en.}, 3(1): 84--85, 1969.
		
		\bibitem[Bla19]{Blanco19}
		G. Blanco.
		\newblock Poles of the complex zeta function of a plane curve.
		\newblock {\em Adv. Math.}, 350: 396--439, 2019.
		
		\bibitem[Bla21]{Blanco31}
		G. Blanco.
		\newblock Yano's conjecture.
		\newblock {\em Invent. Math.}, 226(2): 421--465, 2021.
		
		\bibitem[BLW18]{Milnor_Fibration_Monodromy_3}
		N. Budur, Y. Liu, and B. Wang.
		\newblock The monodromy theorem for compact {K}\"ahler manifolds and smooth
		quasi-projective varieties.
		\newblock {\em Math. Ann.}, 371(3-4): 1069--1086, 2018.
		
		\bibitem[BMS06]{Budur_Mustata_Saito_B_S_for_Ideal}
		N.~Budur, M.~Musta\c{t}\v{a}, and M.~Saito.
		\newblock Bernstein-{S}ato polynomials of arbitrary varieties.
		\newblock {\em Compos. Math.}, 142(3): 779--797, 2006.
		
		\bibitem[BMT11]{Monodromy_Conjecture_for_Hyperplane_Arrangement_Budur_Mustata_Teitler}
		N. Budur, M. Musta\c{t}\u{a}, and Z. Teitler.
		\newblock The monodromy conjecture for hyperplane arrangements.
		\newblock {\em Geom. Dedicata}, 153: 131--137, 2011.
		
		\bibitem[BS10]{Jumping_Coefficient_and_Spectrum_of_a_Hyperplane_Arrangement_Budur_Saito}
		N. Budur and M. Saito.
		\newblock Jumping coefficients and spectrum of a hyperplane arrangement.
		\newblock {\em Math. Ann.}, 347(3): 545--579, 2010.
		
		\bibitem[BSY11]{Local_Zeta_Function_And_b_Function_of_Certain_Hyperplane_Arragement_Budur_Saito_Sergey}
		N. Budur, M. Saito, and S. Yuzvinsky.
		\newblock On the local zeta functions and the {$b$}-functions of certain
		hyperplane arrangements.
		\newblock {\em J. Lond. Math. Soc. (2)}, 84(3): 631--648, 2011.
		\newblock With an appendix by Willem Veys.
		
		\bibitem[BSZ24]{MPVIHA}
		N.~Budur, Q. Shi, and H. Zuo.
		\newblock Motivic principal value integral for hyperplane arrangements.
		\newblock {\em \url{https://arxiv.org/abs/2411.01305}}, 2024.
		
		
		\bibitem[BVWZ21]{Zero_Loci_of_Bernstein-Sato_Ideals}
		N. Budur, R. van~der Veer, L. Wu, and P. Zhou.
		\newblock Zero loci of {B}ernstein-{S}ato ideals.
		\newblock {\em Invent. Math.}, 225(1): 45--72, 2021.
		
		
		\bibitem[DL98]{Motivic_Igusa_Zeta_Function}
		J. Denef and F. Loeser.
		\newblock Motivic {I}gusa zeta functions.
		\newblock {\em J. Algebraic Geom.}, 7(3): 505--537, 1998.
		
		\bibitem[DL99]{Motivic_Integration_on_Arbitrary_Varieties_Denef_Loeser}
		J. Denef and F. Loeser.
		\newblock Germs of arcs on singular algebraic varieties and motivic
		integration.
		\newblock {\em Invent. Math.}, 135(1): 201--232, 1999.
		
		
		\bibitem[DLY24]{DLY24}
		D. Davis and A. C. Lőrincz and R. Yang. 
			\newblock{Archimedean zeta functions, singularities, and Hodge theory}.
			\newblock{\em \url{https://arxiv.org/abs/2412.07849}}, 2024.
		
		
		\bibitem[ELSV04]{Jumping_Coefficients_of_Multiplier_Ideals}
		L. Ein, R. Lazarsfeld, K.~E. Smith, and Dror Varolin.
		\newblock Jumping coefficients of multiplier ideals.
		\newblock {\em Duke Math. J.}, 123(3): 469--506, 2004.
		
		\bibitem[ESV92]{Cohomology_of_Local_System_on_the_Complement_of_Hyperplane_Arrangement}
		H. Esnault, V. Schechtman, and E. Viehweg.
		\newblock Cohomology of local systems on the complement of hyperplanes.
		\newblock {\em Invent. Math.}, 109(3): 557--561, 1992.
		
		\bibitem[Igu74]{Igusa_Zeta_Function_Prologue}
		J. Igusa.
		\newblock Complex powers and asymptotic expansions. {I}. {F}unctions of certain
		types.
		\newblock {\em J. Reine Angew. Math.}, 268/269: 110--130, 1974.
		
		\bibitem[Kas77]{Rationality_of_Roots_of_B-Function_Kashiwara}
		M. Kashiwara.
		\newblock {$B$}-functions and holonomic systems. {R}ationality of roots of
		{$B$}-functions.
		\newblock {\em Invent. Math.}, 38(1): 33--53, 1976/77.
		
		\bibitem[Kas83]{Kashiwara_Malgrange_Theorem_Kashiwara}
		M.~Kashiwara.
		\newblock Vanishing cycle sheaves and holonomic systems of differential
		equations.
		\newblock Algebraic geometry, {Proc}. {Jap}.-{Fr}. {Conf}., {Tokyo} and {Kyoto}
		1982, {Lect}. {Notes} {Math}. 1016: 134-142, 1983.
		
		
		\bibitem[Kol97]{LCT_LARGEST_POLE}
		J. Koll\'ar.
		\newblock Singularities of pairs.
		\newblock In {\em Algebraic geometry---{S}anta {C}ruz 1995}, volume 62, Part 1
		of {\em Proc. Sympos. Pure Math.}, pages 221--287. Amer. Math. Soc.,
		Providence, RI, 1997.
		
		
		\bibitem[Lee24]{lee2024multiplicativethomsebastianibernsteinsatopolynomials}
		J. Lee.
		\newblock Multiplicative Thom-Sebastiani for bernstein-sato polynomials.
		\newblock \url{https://arxiv.org/abs/2402.04512}, 2024.
		
		\bibitem[Lit89]{Lichtin89}
		B. Lichtin.
		\newblock {Poles of {$|f(z,w)|^{2s}$} and roots of the {$b$}-function}.
		\newblock {\em Ark. Mat.}, 27(2): 283--304, 1989.
		
		\bibitem[Loe85]{Quelques}
		F. Loeser.
		\newblock Quelques cons\'equences locales de la th\'eorie de {H}odge.
		\newblock {\em Ann. Inst. Fourier (Grenoble)}, 35(1): 75--92, 1985.
		
		
		
		\bibitem[Mal83]{Kashiwara_Malgrange_Theorem_Malgrange}
		B. Malgrange.
		\newblock Polyn{\^o}mes de {Bernstein}-{Sato} et cohomologie evanescente.
		\newblock Ast{\'e}risque 101-102: 243-267, 1983.
		
		\bibitem[Pop18]{DMBG}
		M.~Popa.
		\newblock {$\mathscr D$}-modules in birational geometry.
		\newblock \url{https://people.math.harvard.edu/~mpopa/notes/DMBG-posted.pdf},
		2018.
		
		\bibitem[PS00]{General_Archimeadean_Zeta}
		D.~H. Phong and J. Sturm.
		\newblock Algebraic estimates, stability of local zeta functions, and uniform
		estimates for distribution functions.
		\newblock {\em Ann. of Math. (2)}, 152(1): 277--329, 2000.
		
		
		\bibitem[Sai94]{microlocal b-fun}
		M. Saito.
		\newblock On microlocal {$b$}-function,
		\newblock{\em Bull. Soc. Math. France}, 122(2): 163--184, 1994.
		
		\bibitem[Sai16]{B_S_for_Centeral_Generic_Hyperplane_Arrangements_Saito}
		M.~Saito.
		\newblock Bernstein-{S}ato polynomials of hyperplane arrangements.
		\newblock {\em Selecta Math. (N.S.)}, 22(4): 2017--2057, 2016.
		
		\bibitem[Sch80]{Milnor_Fibration_Monodromy_2}
		J. Scherk.
		\newblock On the monodromy theorem for isolated hypersurface singularities.
		\newblock {\em Invent. Math.}, 58(3): 289--301, 1980.
		
		
		\bibitem[Ste22]{Mixed_Hodge_Applied_to_Singularity}
		J. Steenbrink.
		\newblock Mixed {H}odge structures applied to singularities.
		\newblock In {\em Handbook of geometry and topology of singularities {III}},
		pages 645--678. Springer, Cham, [2022] \copyright 2022.
		
		\bibitem[SZ24]{shi2024tensorpropertybernsteinsatopolynomial}
		Q. Shi and H. Zuo.
		\newblock On the tensor property of Bernstein-Sato polynomial.
		\newblock \url{https://arxiv.org/abs/2406.04121}, 2024.
		
		\bibitem[Vey24]{Introduction_to_the_Monodromy_Conjecture_Veys_2024}
		W. Veys.
		\newblock Introduction to the monodromy conjecture.
		\newblock \url{https://arxiv.org/pdf/2403.03343}, 2024.
		
		\bibitem[Wal05]{B_S_for_Centeral_Generic_Hyperplane_Arrangements_Walther}
		U. Walther.
		\newblock Bernstein-{Sato} polynomial versus cohomology of the {Milnor} fiber
		for generic hyperplane arrangements.
		\newblock {\em Compos. Math.}, 141(1): 121--145, 2005.
		
		\bibitem[Wal17]{Jacobian module}
		U. Walther.
		\newblock The Jacobian module, the Milnor fiber, and the D-module generated by $f^s$.
		\newblock {\em Invent. Math.}, 207(3): 1239–1287, 2017.
	\end{thebibliography}
\end{document}